\documentclass[11pt,a4paper]{amsart}
\usepackage{amscd}
\theoremstyle{plain} 
\newtheorem{thm}{Theorem}[section]

\newtheorem{lem}[thm]{Lemma}
\newtheorem{prop}[thm]{Proposition}

\theoremstyle{definition}

\newtheorem{rem}[thm]{Remark}

\theoremstyle{remark}

\numberwithin{equation}{section}

\newsavebox{\SmallMathBox}

\def\dd{\partial}

\def\Di{D\kern -.65em /}
\def\Dii{D\kern -.45em /}
\def\di{{\dd}\kern -.55em /}
\def\dii{{\dd}\kern -.40em /}

\def\noi{\noindent}

\def\tand{\mbox{\ \rm  and }}

\def\Del{\Delta}

\def\Cc{{\mathcal C}}
\def\Dd{{\mathcal D}}
\def\Ee{{\mathcal E}}

\def\Gg{{(\ker B)}}
\def\Hh{{\mathcal H}}

\def\Ll{{\mathcal L}}

\def\Nn{{\mathcal N}}

\def\Rr{{\mathcal R}}

\def\={\cong}
\def\>{\supset}
\def\<{\subset}

\def\12{\frac{1}{2}}
\def\2{\Dd}
\def\3{\Nn}
\def\4{\Rr}
\def\6{\cup}
\def\8{\otimes}
\def\0{^{\circ}}

\def\e{\varepsilon}

\def\g{\gamma}
\def\G{\Gamma}

\def\la{\lambda}

\def\s{\sigma}
\def\Si{\Sigma}

\def\z{\zeta}

\def\dom{\mbox{\rm dom\,}}

\def\Id{{\rm Id}}

\def\Si{S\kern -.65em /}
\def\sign{\mbox{\rm sign\,}}

\def\supp{\mbox{\rm supp}}

\def\tr{\mbox{\rm tr\,}}
\def\Tr{\mbox{\rm Tr\,}}

\def\zd{{\det}_{\zeta}}
\begin{document}

\title[Decomposition of the $\z$-determinant and Scattering Theory]
{Adiabatic Decomposition of the $\z$-determinant and Scattering
Theory.}

\author{Jinsung Park}
\address{School of Mathematics\\ Korea Institute for Advanced
Study\\ 207-43\\ Cheongnyangni 2-dong\\ Dongdaemun-gu\\ Seoul
130-722\\ Korea} \email{jinsung@kias.re.kr }

\author{Krzysztof P. Wojciechowski}
\address{Department of Mathematics\\IUPUI (Indiana/Purdue)\\
Indianapolis IN 46202--3216, U.S.A.}
\email{kwojciechowski@math.iupui.edu}

\date{{\em \today. File name:} FinalParkWojMMJ.tex}

\maketitle

\section{Introduction and Statement of the Results}

Let $\Dd : C^{\infty}(M,S) \to C^{\infty}(M,S)$ be a compatible
Dirac operator acting on sections of a Clifford bundle $S$ over a
closed manifold $M$ of dimension $n$. The operator $\Dd$ is a
self-adjoint operator, with discrete spectrum $\{\la_k\}_{k \in
\mathbb{Z}}$. The $\z$-determinant of the Dirac Laplacian $\Dd^2$
is given by the formula
\begin{equation}\label{e:zd}
\zd \Dd^2 = e^{-\z'_{\Dd^2}(0)} \ \, ,
\end{equation}
where $\z_{\Dd^2}(s)$ is defined as follows:
\begin{equation}\label{e:defz}
\z_{\Dd^2}(s) = \sum_{\la_k \ne 0}(\la_k^2)^{-s} = \frac
1{\G(s)}\int_0^{\infty} t^{s-1}[\ \Tr(e^{-t\Dd^2}) - \dim \ker
\Dd\ ]\ dt \ \ .
\end{equation}
This is a holomorphic function of $s$ for $\Re(s)\gg 0$ and has
the meromorphic extension to the complex plane with $s=0$ as a
regular point.

Let us consider a decomposition of $M$ as $M_1 \cup M_2$\,, where
$M_1$ and $M_2$ are compact manifolds with boundaries such that
\begin{equation}\label{e:dec1}
M = M_1 \cup M_2 \, , \ \ M_1 \cap M_2 = Y = {\partial}M_1 =
{\partial}M_2 \, .
\end{equation}
In this paper, we study the adiabatic decomposition of the
$\z$-determinant of $\Dd^2$, which  describes the contributions in
$\zd \Dd^2$ coming from the submanifolds $M_1$ and $M_2$.
Throughout the paper, we assume that the manifold $M$ and the
operator $\Dd$ have product structures in a neighborhood of the
cutting hypersurface $Y$.  Hence, there is a bicollar neighborhood
$N \cong [-1,1]_u \times Y$ of $Y\cong \{0\}\times Y$ in $M$ such
that the Riemannian structure on $M$ and the Hermitian structure
on $S$ are products of the corresponding structures over
$[-1,1]_u$ and $Y$ when restricted to $N$, so that $\Dd$ has the
following form,
\begin{equation}\label{e:pr}
\Dd = G(\partial_u + B)  \quad \text{over} \quad N \ \, .
\end{equation}
Here $u$ denotes the normal variable, $G : S|_Y \to S|_Y$ is a
bundle automorphism, and $B$ is a corresponding Dirac operator on
$Y$. Moreover, $G$ and $B$ do not depend on $u$ and they satisfy
\begin{equation}\label{e:bd}
G^* = - G \ \ , \ \  G^2 = -\Id \ , \ \ B = B^* \ \ \tand \ \ GB =
-BG \ \, .
\end{equation}

To prove the adiabatic decomposition formula of $\zd \Dd^2$, we
follow the original Douglas-Wojciechowski proof of the
decomposition formula for the $\eta$-invariant in \cite{DoWo91}.
However, we face two new problems, not present in the case of the
$\eta$-invariant. First, $\zd \Dd^2$ is a much more
\emph{non-local} invariant than the $\eta$-invariant. This
results, for instance, in the fact that the value of $\zd\Dd^2$
varies with the length of the cylinder. Second, the contribution
of $\zd \Dd^2$ over the cylindrical part is now non-trivial. We
still follow the idea of the paper \cite{DoWo91} and we stretch
our manifold $M$ to separate $M_1$ and $M_2$.  For this, let us
introduce a manifold $M_R$ equal to the manifold $M$ with $N$
replaced by $N_R \cong [-R,R]_u \times Y$.  By assumption of
product structures over $N$, we can extend the bundle $S$ to
$M_R$. Furthermore, using \eqref{e:pr}, we can extend $\Dd$ to the
Dirac operator $\Dd_R$ over $M_R$. Now we decompose $M_R$ by the
hypersurface $\{0\}\times Y$ into two submanifolds $M_{1,R}$,
$M_{2,R}$ and we obtain $\Dd_{1,R}$, $\Dd_{2,R}$ by restricting
$\Dd_R$ to $M_{1,R}$, $M_{2,R}$ respectively.

To formulate the decomposition formula for the $\z$-determinant,
we have to describe the invariant on a manifold with boundary
which enters the picture at this point. The tangential operator
$B$ has discrete spectrum with infinitely many positive and
infinitely many negative eigenvalues. Let $\Pi_>$, $\Pi_<$ denote
the Atiyah-Patodi-Singer (APS) spectral projections onto the
subspaces spanned by the eigensections of $B$ corresponding to the
positive, negative eigenvalues respectively.  We select two
involutions $\sigma_1,\sigma_2$ on kernel of $B$, which satisfy
$G\s_i = - \s_iG$ and define $\pi_i = {\frac{1 - \s_i}{2}}$ the
orthogonal projections onto $-1$ eigenspaces of $\s_i$ . We define
\begin{equation}\label{e:id1}
P_1 = \Pi_< + \pi_1 \ \ , \ \ P_2 = \Pi_> + \pi_2 \ \, ,
\end{equation}
which provides us with the \emph{ideal boundary condition}
introduced by Cheeger in \cite{Ch83}, \cite{Ch87}. The projection
$P_i$ imposes an elliptic boundary condition for $\Dd_{i,R}$ (see
\cite{AtPaSi75}; see \cite{BoWo93} for an exposition of the theory
of elliptic boundary problems for Dirac operators). This means
that the associated operator
$$(\Dd_{i,R})_{P_i}\  =\ \Dd_{i,R} : \dom (\Dd_{i,R})_{P_i} \to L^2(M_{i,R},S)$$
where
$$\dom(\Dd_{i,R})_{P_i} = \{s \in H^1(M_{i,R},S) \mid P_i(s|_Y) =
0\}$$ is a self-adjoint Fredholm operator with $\ker
(\Dd_{i,R})_{P_i} \subset C^{\infty}(M_{i,R},S)$ and discrete
spectrum.

The main concern of this paper is to consider the limit of the
following ratio of the $\z$-determinants,
\begin{equation}\label{e:r1}
{\frac{\zd\Dd_R^2}{\zd(\Dd_{1,R})_{P_1}^2 {\cdot}
\zd(\Dd_{2,R})_{P_2}^2}} \quad \text{as} \quad R \to \infty \ \, ,
\end{equation}
which we call as the \emph{adiabatic decomposition} of the
$\z$-determinant of $\Dd^2$.

The eigenvalues of $\Dd_R$ fall into three different categories as
$R\to\infty$. There are infinitely many large eigenvalues
(\emph{l-values}) bounded away from $0$ and  infinitely many small
eigenvalues (\emph{s-values}) of the size $O({R}^{-1})$ . Besides
these, there are finitely many eigenvalues which decay
exponentially with $R$ (\emph{e-values}). The number $h_M$ of
\emph{e-values} is given by
\begin{equation}\label{e:e1}
h_M =  \dim \ker_{L^2}\Dd_{1,\infty} + \dim
\ker_{L^2}\Dd_{2,\infty} + \dim  L_1 \cap L_2 \ \, ,
\end{equation}
where $\Dd_{i,\infty}$ is the operator defined from $\Dd$ in a
natural way over the manifold $M_{i,\infty}$ , which is equal to
$M_i$ with the half infinite cylinder $[0,\infty) \times Y$ or
$(-\infty,0]\times Y$ attached. More precisely, the operator
$\Dd_i=\Dd|_{M_i}$ extends in a natural way to the manifold
$M_{i,\infty}$ . It has a unique closed self-adjoint extension in
$L^2(M_{i,\infty},S)$ , which we denote by $\Dd_{i,\infty}$ . The
subspaces $L_i \subset \ker B$ are the spaces of limiting values
of extended $L^2$-solutions of $\Dd_{i,\infty}$. The decomposition
of the eigenvalues of the operator $\Dd_R$ into different classes
was discussed by Cappell, Lee and Miller (see \cite{CLM1}). The
corresponding analysis for the operator $(\Dd_{i,R})_{P_i}$ was
provided by M${\ddot {\rm u}}$ller (see \cite{Mu94}). The spectrum
of the operators $(\Dd_{i,R})_{P_i}$ splits in the same way as the
spectrum of $\Dd_R$ . The only difference is that the operators
$(\Dd_{i,R})_{P_i}$ do not have nonzero \emph{e-values} and the
dimension of the space of the solutions of $(\Dd_{i,R})_{P_i}$ is
equal to
\begin{equation}\label{e:e2}
h_i = \dim \ker  (\Dd_{i,R})_{P_i} = \dim \ker_{L^2}\Dd_{i,\infty}
+ \dim  L_i \cap \ker(\s_i - 1) \ \, .
\end{equation}

In the adiabatic limit process, the different types of eigenvalues
make their contribution at different time intervals of the
integral representation of $\z_{\Dd^2}(s)$ in \eqref{e:defz}. The
contribution made by \emph{l-values} comes from the time interval
$[0, R^{2-\e}]$ , where $\e$ is a sufficiently small positive
number, and we fix $\e$ from now on. More precisely, it is not
difficult to show that the \emph{l-values} contribution to the
adiabatic limit of (\ref{e:r1}) from the time interval $[R^{2-\e},
\infty]$ disappears as $R \to \infty$ (see Section \ref{two}). The
contribution made by \emph{l-values} was discussed in
\cite{JKPW2}. To be more precise, in \cite{JKPW2} we discussed the
case of the operator $\Dd$ , such that $\Dd_{i,\infty}$ and $B$
have trivial kernels. These conditions imply that there are no
\emph{e-values} and \emph{s-values}. This allows us to reduce the
computation of the quotient in (\ref{e:r1}) to the corresponding
quotient on the cylinder, hence one can show that the limit of
(\ref{e:r1}) as $R \to \infty$ is equal to $2^{-\z_{B^2}(0)}$ .
Actually, even in the presence of \emph{e-values} and
\emph{s-values}, we are able to show that in the adiabatic limit
the contribution of \emph{l-values} comes only from the time
interval $[0,R^{2-\e}]$ so that we can reduce to the cylinder as
in \cite{JKPW2}. The method we use to prove this combines
{Duhamel's principle} and {Finite propagation speed property of
the wave operators}. Details are presented in Section \ref{two}.

The \emph{s-values} contribution comes from the time interval
$[R^{2-\e}, R^{2+\e}]$. The computation of the \emph{s-values}
contribution is the main achievement of this paper. We follow
M$\ddot {\rm u}$ller (see \cite{Mu94}) and use the
\emph{Scattering theory} to get a description of the
\emph{s-values}. The operators $\Dd_{i,\infty}$ on $M_{i,\infty}$
determine scattering matrices $C_{i}(\la)$ . It turned out that
the matrix $C_{12} = C_1(0)\circ C_2(0)$ on $\ker B\cap \ker(G+i)$
determines the contribution of \emph{s-values} of the operator
$\Dd_R$ in the adiabatic limit. Similarly  the finite-dimensional
unitary matrix $S_{\s_i}$ on $\ker(\sigma_i+1)$, which is defined
by the scattering matrix $C_i(0)$ and the involution $\sigma_i$,
determines the contribution of \emph{s-values} of the operators
$(\Dd_{i,R})_{P_i}$. The exact correspondence is stated in Section
\ref{three}.

Finally we have to discuss \emph{e-values} of $\Dd_R$ . The number
of \emph{e-values} is equal to $h_M$ which, as remarked above, is
constant. On the other hand, the set of zero eigenvalues of
$\Dd_R$, which is a subset of \emph{e-values} by definition, is
very unstable with respect to $R$ . Hence, without making
additional assumptions we are not able to control the adiabatic
limit of the determinant of $\Dd_R^2$ due to the finite number of
nonzero \emph{e-values}. Hence, we assume that all the
\emph{e-values} are zero eigenvalues in order to avoid the
technical difficulty of the nonzero \emph{e-values}. One of the
important examples of such situation is the case of the operator
$$
d_{\rho}+d_{\rho}^* \ : \ \oplus_{i=0}^{n} \Omega^i(M, V_{\rho}) \
\rightarrow \ \oplus_{i=0}^{n} \Omega^i(M, V_{\rho})
$$
where $V_{\rho}$ denotes the flat vector bundle defined by the
unitary representation $\rho$ of $\pi_1(M)$ (see Proposition
\ref{p:zero}). For the operator $L : W \to W$ acting on a finite
dimensional vector space $W$, we denote by ${\det}^* L$ the
determinant of the operator $L$ restricted to the subspace $(\ker
L)^{\perp}$. Now we are ready to formulate the main result of the
paper.

\begin{thm}\label{t:apssplit} When all the \emph{e-values} of $\Dd_R$ are zero eigenvalues,
the following formula holds:
\begin{multline}\label{e:main}
 \lim_{R \to \infty} R^{-2h}\cdot \frac{{\zd\Dd_R^2}}
{\zd(\Dd_{1,R})_{P_1}^2{\cdot}\zd{(\Dd_{2,R})_{P_2}^2}}  \\ \ =
2^{-\z_{B^2}(0)-h_Y+2h_M}{\cdot} {\det}^* \biggl(\frac{2\, \Id
-C_{12}-C_{12}^{-1}}{4}\biggr)
  \cdot \prod^2_{i=1}{\det}^* \biggl(\frac{2\, \Id-
S_{\sigma_i}-S_{\sigma_i}^{-1}}{4}\biggr)^{-1} \ \,
\end{multline}
where $h=h_M-h_1-h_2$ and $h_Y = \dim \ker B $ .
\end{thm}

\begin{rem}
In \cite{HMaM95} and \cite{Ha98}, the reduced normal operators
corresponding to our operators $C_{12}$, $S_{\sigma_i}$ were
introduced in the framework of {\it b-calculus} and used in the
analysis of {\it s-values} for the analytic surgery of the
$\eta$-invariant and analytic torsion.
\end{rem}

To prove Theorem \ref{t:apssplit}, we consider the following
relative $\z$-function and its derivative at $s=0$,
\begin{align*}\label{e:rel1}
\frac{1}{\Gamma(s)}\int^{\infty}_{0}t^{s-1} [\, \Tr(e^{-t\Dd_R^2}
- e^{-t(\Dd_{1,R})_{P_1}^2} - e^{-t(\Dd_{2,R})_{P_2}^2})- h\, ] \
dt \ \, ,
\end{align*}
which we decompose  into two parts,
$$\z_{s}^R(s)=\frac{1}{\Gamma(s)}\int^{R^{2-\e}}_0 (\cdot)\ dt \ \ ,
\ \ \z_l^R(s)=\frac{1}{\Gamma(s)}\int^{\infty}_{R^{2-\e}} (\cdot)\
dt \ \ ,$$ where $\e$ is the fixed sufficiently small positive
number. The derivatives of $\z_s^R(s)$ and $\z_l^R(s)$ at $s=0$
give the small and large time contribution in (\ref{e:main}).

In Section \ref{two} we deal with the small time contribution and
prove that this equal $2^{-\z_{B^2}(0)}$, which gives  the first
factor on the right side of (\ref{e:main}). In Section \ref{three}
we explain some basic description of the small eigenvalues. We
follow \cite{Mu94} and use scattering theory in order to get a
description of the \emph{s-values} of $\Dd_R$ and
$(\Dd_{i,R})_{P_i}$ , which allows us to make a comparison of
\emph{s-values} of those operators with the eigenvalues of certain
model operators over $\mathbb{S}^1$ . This is the central part of
this paper. In Section \ref{four} we use the results of Section
\ref{three} to show that, in the adiabatic limit, the large time
contribution to the quotient (\ref{e:r1}) is equal to
$$
 2^{-h_Y+2h_M}{\cdot} {\det}^*
\biggl(\frac{2\mathrm{Id}-C_{12}-C_{12}^{-1}}{4}\biggr)
  \cdot \prod^2_{i=1}{\det}^* \biggl(\frac{2\mathrm{Id}-
S_{\sigma_i}-S_{\sigma_i}^{-1}}{4}\biggr)^{-1} \ \ .
$$
This is the second factor on the right side of (\ref{e:main}). The
zero eigenvalues make their presence via the factor $R^{-2h}$ on
the left side of (\ref{e:main}).

In Section \ref{five} we review the decomposition formula for the
$\eta$-invariant and offer a new proof based on the method
developed in order to prove Theorem \ref{t:apssplit}. This proof
is more complicated than other proofs presented in \cite{Bu95},
\cite{Dfr94}, \cite{HMaM95}, \cite{Mu96}, \cite{KPW95},
\cite{BL98}, \cite{KL01}, \cite{LP1}. However, it is a nice
illustration of the differences we encounter when we deal with the
$\z$-determinant instead of the $\eta$-invariant.

{\bf Acknowledgments.} The authors want to express their gratitude
to Werner M${\ddot {\rm u}}$ller for his helpful comments on the
content of this paper and for the more general remarks on the
theory of the $\z$-determinant.  A part of this work was done
during the first author's stay at ICTP and MPI. He wishes to
express his gratitude to ICTP and MPI for their financial support
and hospitality.

\section{Small Time Contribution}
\label{two}

In this section we determine the small time contribution, which is
done in two steps. First, we use Duhamel's principle and Finite
propagation speed property of the wave operator to show that we
can reduce the problem to computations on the cylinder. Then, we
perform the explicit calculations on the cylinder. Both parts are
fairly standard. The cylinder contribution has been recently
computed in \cite{JKPW2}. Therefore, we only discuss the reduction
scheme and refer to \cite{JKPW2} for the explicit computation on
the cylinder.

Let $\Ee_R(t;x,y)$ denote the kernel of the operator
$e^{-t\Dd_R^2}$ . We introduce the specific parametrix for
$\Ee_R(t;x,y)$ , which fits our main purpose to \emph{localize}
the contribution coming from the cylinder $[-R,R]_u \times Y$ and
the interior of $M_R$ . In  fact, the interesting point here is
that we use $\Ee_R(t;x,y)$ to construct this parametrix. Let
$\Ee_{{c}}(t;x,y)$ denote the kernel of the operator
$e^{-t(-\partial_u^2 + B^2)}$ on the infinite cylinder $\mathbb{R}
\times Y$ . We introduce a smooth, increasing function $\rho(a,b)
: [0, \infty) \to [0,1]$ equal to $0$ for $0 \le u \le a$ and
equal to $1$ for $b \le u$ . We use $\rho(a,b)(u)$ to define
\begin{align*}
 \phi_{1,R} = 1 - \rho(\frac{5}{7}R , \frac{6}{7}R) \ \ , \ \
 \psi_{1,R}
= 1 - \psi_{2,R} \ \ , \\ \phi_{2,R} = \rho(\frac{1}{7}R ,
\frac{2}{7}R) \ \ , \ \ \psi_{2,R} = \rho(\frac{3}{7}R ,
\frac{4}{7}R) \ \  .
\end{align*}
We extend these functions to symmetric functions on the whole real
line. These functions are constant outside the interval $[-R ,
R]_u$ and we use them to define the corresponding functions on a
manifold $M_R$, which are denoted by the same notations. Now, we
define $Q_R(t;x,y)$ a \emph{parametrix} for the kernel
${\Ee}_R(t;x,y)$ by
\begin{equation}\label{e:p1}
Q_R(t;x,y) = \phi_{1,R}(x){\Ee}_{c}(t;x,y)\psi_{1,R}(y) +
\phi_{2,R}(x){\Ee}_R(t;x,y)\psi_{2,R}(y) \ \, .
\end{equation}
It follows from the Duhamel's principle that
\begin{equation}\label{e:p2}
{\Ee}_R(t;x,y) = Q_R(t;x,y) + ({\Ee}_R*\Cc_R)(t;x,y) \ \, ,
\end{equation}
where ${\Ee}_R * \Cc_R$ is the convolution given by
$$({\Ee}_R * \Cc_R)(t;x,y) = \int_0^tds\int_{M_R}\, dz\,
\Ee_R(s;x,z)\, \Cc_R(t-s;z,y)  \ \ ,$$ and the error term
$\Cc_R(t;x,y)$ is given by the formula
\begin{align*}
\Cc_R(t;x,y) =&  - \partial_u^2\phi_{1,R}(x)\, {\Ee}_{c}(t;x,y)\,
\psi_{1,R}(y) - {\partial}_u\,\phi_{1,R}(x)\,
{\partial}_u{\Ee}_{c}(t;x,y)\,\psi_{1,R}(y)\\
& -{\partial}_u^2\phi_{2,R}(x)\, {\Ee}_R(t;x,y)\,\psi_{2,R}(y) -
{\partial}_u\phi_{2,R}(x)\, {\partial}_u{\Ee}_R(t;x,y)\,
\psi_{2,R}(y) \  .
\end{align*}
The following elementary lemma follows from the construction of
$Q_R(t;x,y)$,

\begin{lem}\label{l:err1}For a fixed $y$,
the support of the error term $\Cc_R(t;x,y)$ as a function of $x$
is a subset of $\big([-\frac{6}{7}R , -\frac{1}{7}R]_u\cup
[\frac17 R, \frac{6}{7}R]_u\big) \times Y$ . Moreover it is equal
to $0$ if the distance between $x$ and $y$ is smaller than $\frac
R7$ .
\end{lem}

 Now, following Cheeger, Gromov and Taylor (see \cite{CGT82}; see
also Section 3 of \cite{Bu95}), we use the Finite propagation
speed property for the wave operator. The technique introduced in
\cite{CGT82} allows us to compare the heat kernel of the operator
$\Dd_R^2$ over $M_R$ with the heat kernel of the operator
$-\partial_u^2 + B^2$ on the cylinder $\mathbb{R} \times Y$ . We
describe the case we need in our work. Let $X_1$ and $X_2$ be
Riemannian manifolds of dimension $n$ and $S_i$ be spinors bundle
with Dirac operators $\Dd_i$ over $X_i$. Assume that there exists
a decomposition $X_i = K_i \cup U_i$ , where $U_i$ is an open
subset of $X_i$ . Moreover, we assume that there exists an
isometry $h : U_1 \to U_2$ covered by the unitary bundle
isomorphism $\Phi_h : S_1|_{U_1} \to S_2|_{U_2}$ , which
intertwines Dirac operators $\Dd_1|_{U_1}$ and $\Dd_2|_{U_2}$ . We
identify
$$U \cong U_1 \cong U_2 \ \ ,$$
so that $X_1$ and $X_2$ have a common open subset $U$ . Let
$\Ee_i(t;x,y)$ denote the kernel of the operator $e^{-t\Dd_i^2}$ .
Then we have the following estimate on the difference of the heat
kernels on $U$ as in Lemma 3.6 in \cite{Bu95},

\begin{prop}\label{p:fps0}
For $x, y \in U$ and $t>0$, there exist positive constants $c_1,
c_2$ such that
\begin{equation}\label{e:fps}
\| \partial_u^j \Ee_1(t;x,y) -\partial_u^j \Ee_2(t;x,y) \| \le
c_1e^{-c_2{\frac {r^2}t}} \ \,
\end{equation}
where $j=0,1$ and $r = \min(d(x,K_1), d(y,K_1))$.
\end{prop}

In our situation, $X_1 = M_R$, $X_2=\mathbb{R} \times Y$ and
$U=[-R,R]_u\times Y$. Note that the heat kernel $\Ee_{c}(t;x,y)$
over $\mathbb{R}\times Y$ satisfies the standard estimate. More
precisely, for $t>0$ we have
\begin{equation}\label{e:st}
\|\partial_u^j \Ee_{c}(t;(u,w),(v,z))\| \le c_1|u-v|^j t^{- \frac
n2-j} e^{-c_3{\frac {(u-v)^2}{t}}} \ \, ,
\end{equation}
where $j=0,1$, $u,v \in \mathbb{R}$ and $w,z \in Y$ . This follows
from the corresponding estimate for the heat kernel of $B^2$ over
the closed manifold $Y$ (see Proposition 4.1 in \cite{JKPW4}) and
explicit form of the heat kernel of $-\partial_u^2$ over
$\mathbb{R}$. We are going to use (\ref{e:fps}) and (\ref{e:st})
in the following proposition.
\bigskip

\begin{prop}\label{p:fps}
There exist constants $c_1,c_2 >0$ such that for any $t$ with
$0<t<R^{2-\e}$ and $((u,w),(v,z)) \in \supp\, \Cc_R(t; \cdot ,
\cdot)$,
\begin{align}\label{e:fps2}
\|\Ee_R(t;(u,w),(v,z))\| \le c_1e^{-c_2{\frac {R^2}t}} \ ,\qquad
\|\Cc_R(t;(u,w),(v,z))\| \le c_1e^{-c_2{\frac {R^2}t}} \ \, .
\end{align}
\end{prop}

\begin{proof}
For $j=0,1$, we have
\begin{multline*}
\|\partial_u^j\Ee_R(t;(u,w),(v,z))\|\  \le  \
  \|\partial^j_u\Ee_{c}(t;(u,w),(v,z))\| \\
+ \|\partial^j_u\Ee_R(t;(u,w),(v,z)) -\partial^j_u
\Ee_{c}(t;(u,w),(v,z))\| \ \ .
\end{multline*}
By  Lemma \ref{l:err1} and \eqref{e:fps}, \eqref{e:st}, there
exist some constants $c_1, c_2>0$ such that for $(u,w),(v,z)  \in
\supp\, \Cc_R(t; \cdot , \cdot)$, both summands on the right side
satisfy the desired estimate.  This estimate for $j=0$ ($j=1$)
implies the first (second) estimate in \eqref{e:fps2}.
\end{proof}

Now we are ready to prove the following technical result,

\begin{prop}\label{p:fps2}
\begin{equation}\label{e:err2}
\lim_{R \to \infty} \frac{d}{ds}\Big|_{s=0}
\frac{1}{\Gamma(s)}\int_0^{R^{2-\e}} t^{s-1} {dt}\int_{M_R}\tr
({\Ee}_R*\Cc_R)(t;x,x)\ dx = 0 \ \, .
\end{equation}
\end{prop}

\begin{proof}
By Lemma \ref{l:err1} and Proposition \ref{p:fps},
\begin{align*}
& \ |\, \tr  ({\Ee}_R*\Cc_R)(t;x,x)\, | \ \le \
\|({\Ee}_R*\Cc_R)(t;x,x)\|\\
\le & \ \int_0^tds\int_{[-\frac{6}{7}R , \frac{6}{7}R]_u \times Y}
\|{\Ee}_R(s;x,z)\, \Cc_R(t-s;z,x)\|\ dz \\
 \le & \ c_1^2{\cdot}\int_0^tds
\int_{[-\frac{6}{7}R , \frac{6}{7}R]_u \times Y}
e^{-c_2{\frac{R^2}{s}}}e^{-c_2{\frac{R^2}{t-s}}}\ dz \\
\le & \ c_3R{\cdot}\int_0^t e^{-c_2{\frac{tR^2}{s(t-s)}}}ds \ \le
\ c_3R{\cdot}\int_0^{\frac t2}e^{-2c_2{\frac{R^2}{s}}}ds\ \le \
c_3R\,{\frac t2}\, e^{-4c_2{\frac{R^2}{t}}} \ ,
\end{align*}
where the last estimate is a consequence of the elementary
inequality
$$\int_0^te^{-\frac cs}ds \le t\, e^{-\frac ct} \ \ .$$
Hence we have proved
\begin{equation}\label{e;fps3}
|\,\tr \ ({\Ee}_R*\Cc_R)(t;x,x)\,| \le c_4 R\, t\, e^{-c_5\frac
{R^2}t} \ \, .
\end{equation}
This allows us to estimate as follows
\begin{multline*}
 \Big| \frac{1}{\Gamma(s)}\int_0^{R^{2-\e}} t^{s-1}{dt} \int_{M_R}
 \tr ({\Ee}_R*\Cc_R)(t;x,x) \ dx \Big| \\ \le \ c_6R^2{\cdot}\,
\Big|\frac{1}{\Gamma(s)}\Big|\int_0^{R^{2-\e}}|t^s| e^{-c_5\frac
{R^2}t}\ dt \ \ .
\end{multline*}
As $R\to\infty$, the function of $s$ on the right side uniformly
converges to zero for $s$ in any compact set in $\mathbb{C}$.
Hence, the derivative at $s=0$ of the meromorphic function on the
left side converges to zero as $R\to\infty$. This completes the
proof.
\end{proof}

The corresponding result for the operator $(\Dd_{i,R})^2_{P_i}$ to
Proposition \ref{p:fps2} can be carried out in exactly the same
manner. First, as for $\Dd_R$ over $M_R$, we can construct the
parametrices for the heat kernels of $e^{-t(\Dd_{i,R})^2_{P_i}}$
using $\Ee_R(t;x,y)$ and the heat kernels of
$\big(G(\partial_u+B)\big)^2_{P_i}$ over $[0,\infty)_u\times Y$ or
$(-\infty, 0]_u\times Y$. Second, one can obtain the corresponding
estimate to \eqref{e:fps} using the explicit form of the heat
kernel of $(G(\partial_u+B))^2_{P_i}$. Third, one can also have
the corresponding estimate to Proposition \ref{p:fps} for
$(\Dd_{i,R})^2_{P_i}$ since the similar estimate as in Proposition
\ref{p:fps0} holds over the support of the error terms. (see Lemma
3.6 in \cite{Bu95}). All these imply that the similar estimate to
Proposition \ref{p:fps2} holds for $(\Dd_{i,R})^2_{P_i}$ . Now we
are ready to prove the following main result of this section,

\begin{prop}\label{p:s}
\begin{equation}\label{e:s}
\lim_{R \to \infty} \Bigr({\frac{d}{ds}}\Big|_{s=0}\z_{s}^R(s) +
h(\g + (2-\e){\cdot}\log R) \Bigr) = \z_{B^2}(0)\cdot \log 2 \ \,
.
\end{equation}
\end{prop}

\begin{proof}

 First we
observe
\begin{align}\label{contr1}
 \frac{d}{ds}\Bigr|_{s=0}  \biggl( \frac{h}{\Gamma(s)}
\int^{R^{2-\e}}_{0}t^{s-1} \ dt\biggr) =  h(\gamma + {(2-\e)}\log
R) \ \ .
\end{align}
Hence we need to compute the limit as $R\to\infty$
 of the following remaining part of ${{\z}_s^R}(s)$,
\begin{align*}
\frac{d}{ds}\Big|_{s=0}\frac{1}{\Gamma(s)}\int^{R^{2-\e}}_{0}t^{s-1}
\Tr(e^{-t\Dd_R^2} - e^{-t(\Dd_{1,R})_{P_1}^2} -
e^{-t(\Dd_{2,R})_{P_2}^2})\ dt.
\end{align*}
By Proposition \ref{p:fps2} and corresponding results for
$(\Dd_{i,R})_{P_i}^2$, it is sufficient to consider the limit as
$R\to\infty$ of
\begin{align*}
\frac{d}{ds}\Big|_{s=0} \frac{1}{\Gamma(s)}\int^{R^{2-\e}}_0
t^{s-1} dt \int_{M_R} \tr \big( Q_R(t;x,x) -  Q_{1,R} (t;x,x) -
Q_{2,R}(t;x,x)\big) dx
\end{align*}
where $Q_{i,R}(t;x,y)$ denotes the parametrix for
$e^{-t(\Dd_{i,R})_{P_i}^2}$. Now, the interior contributions to
the different parametrices, all determined by the kernel
$\Ee_R(t;x,y)$ , cancel out and we are left only with the cylinder
contribution. Hence we have to deal with the limit as $R\to\infty$
of
\begin{multline}\label{e:cylinderpart}
\frac{d}{ds}\Big|_{s=0} \frac{1}{\Gamma(s)}\int^{R^{2-\e}}_0
t^{s-1} dt \int_{M_R} \tr \big( \psi_{1,R}\Ee_c(t;x,x) \\-
\psi_{1,R} \Ee_{c,1} (t;x,x) - \psi_{1,R}\Ee_{c,2}(t;x,x)\big) dx
\end{multline}
where $\Ee_{c,i}(t;x,y)$ denotes the heat kernel of
$(G(\partial_u+B))_{P_i}^2$ over the half cylinder. We repeat
computations in Section 2 of \cite{JKPW2} where we assumed the
conditions that $B$ is invertible and $Y$ is even dimensional. But
we can easily derive the same formula following Section 2 of
\cite{JKPW2} without these assumptions. So we can show that for
$s$ in a compact subset of $\mathbb{C}$ the integral part in
\eqref{e:cylinderpart} uniformly converges to the following
function as $R\to\infty$ ,
\begin{align*}
2\,\Bigl({\frac{\G(s)}{4}}-{\frac{\G(s +
{\frac1{2}})}{4s\sqrt{\pi}}}\Bigr){\cdot}\z_{B^2}(s) \ \, .
\end{align*}
Hence, we obtain
\begin{align}\label{contr2}
&\lim_{R\to\infty}\frac{d}{ds}\Big|_{s=0}\frac{1}{\Gamma(s)}\int^{R^{2-\e}}_{0}t^{s-1}
\Tr(e^{-t\Dd_R^2} - e^{-t(\Dd_{1,R})_{P_1}^2} -
e^{-t(\Dd_{2,R})_{P_2}^2})\ dt\\
&\qquad =\  \frac{d}{ds}\Big|_{s=0}\frac{2}{\Gamma(s)}\,\Bigl(
{\frac{\G(s)}{4}}-{\frac{\G(s + {\frac1{2}})}{4s\sqrt{\pi}}}
\Bigr){\cdot}\z_{B^2}(s) =  \z_{B^2}(0)\cdot \log 2\ \, . \notag
\end{align}
Combining \eqref{contr1} and \eqref{contr2} completes the proof.

\end{proof}

\section{small eigenvalues and scattering matrices}
\label{three}

In this section we investigate the relation between the
\emph{s-values} of the operators $\Dd_R$, $(\Dd_{i,R})_{P_i}$ and
the scattering matrices $C_i(\lambda)$ determined by the operators
$\Dd_{i,\infty}$ on $M_{i,\infty}$ for $i=1,2$. We refer to
Section 4 and Section 8 in \cite{Mu94} for a more detailed
exposition of the elements of \emph{Scattering theory} that we use
in this paper.

Let us recall that $M_R$ has the cylindrical part
$N_R=[-R,R]_u\times Y$. Hence $M_{1,R}, M_{2,R}$ have the
cylindrical part $[-R,0]_u\times Y$, $[0,R]_u\times Y$
respectively. But, in order to consider $M_{i,R}$ as a submanifold
of $M_{i,\infty}$ which is obtained by  attaching
$[0,\infty)_v\times Y$ or $(-\infty,0]_v\times Y$ to $M_{i}$, we
change the variable by $v=u+R$ or $v=u-R$ so that the cylindrical
part of $M_{i,R}$ is given by $[0,R]_v\times Y$ or $[-R,0]_v\times
Y$. Throughout this section, we will use this convention when it
is needed.

 For any $\psi \in
\ker B$ and $\la \in \mathbb{C} - (-\infty, -\mu_1] \cup
[\mu_1,+\infty)$ where $\mu_1$ denotes the lowest positive
eigenvalue of the tangential operator $B$,  there exists a
generalized eigensection $E(\psi,\lambda)$ of $\Dd_{1,\infty}$
over $M_{1,\infty}$ determined by the couple $(\psi,\la)$ (see
Section 4 in \cite{Mu94} ) in the following sense,
$$
\Dd_{1,\infty} E(\psi,\la) = \la\, E(\psi,\lambda) \ \, .$$ The
section $E(\psi,\la)$ has the following form over
$[0,\infty)_v\times Y$,
\begin{equation}\label{E-exp}
E(\psi,\lambda)= e^{-i\lambda v}(\psi-iG\psi) + e^{i\lambda
v}C_1(\lambda)(\psi-iG\psi) + \theta(\psi,\lambda)
\end{equation}
 where $\theta$
is a square integrable section such that, for each $v$ ,
$\theta(\psi,\lambda,(v,\cdot))$ is orthogonal to $\ker B$. The
operator $C_1(\lambda) : \ker B \to \ker B$ is regular and unitary
for $|\lambda|< \mu_1$ and equals the \emph{Scattering matrix}
such that
$$
C_1(\lambda)C_1(-\lambda)=\Id \  , \ \
C_1(\lambda)G=-GC_1(\lambda) \ \, ,
$$
which imply
$$
C_1(0)^2=\Id \ , \ \ C_1(0)G=-GC_1(0) \ \ .$$ Therefore $C_1(0)$
gives a distinguished unitary involution of $\ker B$. In fact, the
space of the limiting values of the extended $L^2$-solutions of
$\Dd_{1,\infty}$, $L_1 \subset \ker B$ is equal to the
$(+1)$-eigenspace of $C_1(0)$, that is, $L_1=\ker(C_1(0)- 1)$. The
following proposition is a basic tool to deal with $E(\psi,\la)$,

\begin{prop}\label{Maass-Selberg}\emph{(Maass-Selberg)} The following equality holds,
\begin{multline*}
\langle E(\phi,\la), E(\psi,\la) \rangle_{M_{1,R}}\\ = 4R \langle
\phi, \psi \rangle_Y - i \langle C_1(-\la) C_1'(\la)(\phi-iG\phi),
\psi-iG\psi \rangle_Y +O(e^{-cR}) \ \
\end{multline*}
where $\phi,\psi\in \ker B$.
\end{prop}

\begin{proof}
By Green's formula, we have
\begin{align}\label{e:MS}
& \ h \langle E(\phi,\la+h), E(\psi,\la) \rangle_{M_{1,R}} \\=&\
\langle \Dd_{1,R} E(\phi,\la+h), E(\psi,\la) \rangle_{M_{1,R}} -
\langle E(\phi,\la+h), \Dd_{1,R} E(\psi,\la)
\rangle_{M_{1,R}} \notag \\
= &\ \langle G E(\phi,\la+h)|_{\partial(M_{1,R})},
E(\psi,\la)|_{\partial(M_{1,R})} \rangle_{\partial(M_{1,R})} \ \,
. \notag
\end{align}
Using \eqref{E-exp}, the last line has the following form,
\begin{align*}
& \ \ i\, e^{-i h R} \langle \phi-iG\phi, \psi-iG\psi \rangle_Y \\
& - i\, e^{i h R}
\langle C_1(\la+h)(\phi-iG\phi), C_1(\la)(\psi-iG\psi) \rangle_Y + O(e^{-cR})\\
=& \ \ i\, e^{-i h R} \langle \phi-iG\phi, \psi-iG\psi \rangle_Y
  - i\, e^{i h R} \langle \phi-iG\phi, \psi-iG\psi \rangle_Y\\
&+ i\, e^{i h R} \langle C_1(\la)(\phi-iG\phi), C_1(\la)
(\psi-iG\psi) \rangle_Y\\
& - i\, e^{i h R} \langle C_1(\la+h)(\phi-iG\phi),
C_1(\la)(\psi-iG\psi) \rangle_Y + O(e^{-cR}) \ \, .
\end{align*}
Now, dividing the right side by $h$ and taking the limit $h\to 0$,
we obtain
\begin{align*}
&\, 2R \langle \phi-iG\phi,\psi-iG\psi \rangle_Y - i \langle
C_1'(\la)( \phi-iG\phi),C_1(\la) (\psi-iG\psi) \rangle_Y +
O(e^{-cR})\\
=&\, 4R \langle \phi,\psi \rangle_Y - i \langle C_1(-\la)
C_1'(\la)( \phi-iG\phi), (\psi-iG\psi) \rangle_Y + O(e^{-cR})\ \,
.
\end{align*}
Comparing this with \eqref{e:MS} (divided by $h$) completes the
proof.

\end{proof}

Now we shall analyze the \emph{s-values} of $\Dd_R$ over $M_R$.
Let us consider a \emph{s-value} $\la=\la(R)$ of $\Dd_R$ such that
$$|\la(R)|\le R^{-\kappa} \quad \text{for sufficiently large $R$} $$
where $\kappa$ is a fixed constant with $0<\kappa<1$. Let $\Psi_R$
denote a \emph{normalized} eigensection of $\Dd_R$ corresponding
to \emph{s-value} $\la$, that is,
$$
\Dd_R\Psi_R =\la\Psi_R   \ \, , \ \ \|\Psi_R\| = 1\ \, .$$ Over
the cylindrical part $[-R,R]_u\times Y$ in $M_R$, the eigensection
$\Psi_R$ corresponding to \emph{s-value} $\la$ of $\Dd_R$ has the
following form,
\begin{equation}\label{e:psi}
\Psi_R = e^{-i\la u}\psi_1 + e^{i\la u}\psi_2 +\hat{\Psi}_R
\end{equation}
where $\psi_1\in \ker B\cap \ker(G-i)$, $\psi_2\in \ker B \cap
\ker(G+i)$ and $\hat{\Psi}_R$ is orthogonal to $\ker B$.
\begin{lem} \label{l:cut-est}
We have the following estimates
\begin{align*}
||\hat{\Psi}_R|_{\{u\}\times Y}||_Y \leq  c_1 e^{-c_2R} \qquad
\text{for}\ \ -\frac34 R \le u \le \frac34 R
\end{align*}
where $c_1,c_2$ are positive constants independent of $R$ .
\end{lem}

The proof of this lemma is same as the one of Lemma 2.1 in
\cite{KPW94}. Now we can prove

\begin{prop}\label{p:nts}
The zero eigenmode $e^{-i\lambda u}\psi_1+e^{i\lambda u}\psi_2$ of
the eigensection $\Psi_R$ of \emph{s-value} $\la(R)$ of $\Dd_R$ is
non-trivial.
\end{prop}

\begin{proof}We follow the proof of Theorem 2.2 in
\cite{KPW94}, so we assume that the zero eigenmode of $\Psi_R$ is
trivial, which will contradict to the fact $\la(R)$ is a
\emph{s-value}. Throughout the proof, we regard $M_{1,R}$ as a
submanifold of $M_{1,\infty}$ using the change of variable
$v=u+R$. We define a section $\Phi_R$ on $M_{1,\infty}$ by
\[
\Phi_R = \begin{cases}& h_R(x)\Psi_R(x) \quad \text{for} \quad
x\in
M_{1,R} \\
& 0 \quad\qquad \text{for}\quad x\in M_{1,\infty}\setminus M_{1,R}
\end{cases}
\]
where $h_R$ is a smooth function on $M_{1,\infty}$, equal to $1$
for $x\in M_1\cup[0, \frac R2]_v\times Y$ and equal to $0$ for
$x\in [\frac {3}4 R,\infty)_v\times Y$ with $|\frac{\partial^j
h}{\partial v^j}| \le CR^{-j}$ for a constant $C>0$. Let
$H^1(M_{1,\infty},S)$ denote the first Sobolev  space. For any
$a\ge 0$, we introduce a closed subspace of $H^1(M_{1,\infty},S)$
by
\begin{multline*}
H^1_a(M_{1,\infty},S)\\=\{ \Phi\in H^1(M_{1,\infty},S) \ | \
\langle \Phi(v,\cdot), \phi_k \rangle =0 \quad \text{for}\quad
v\ge a, \ k=1,\ldots, h_Y\}
\end{multline*}
where $\phi_1,\ldots, \phi_{h_Y}$ denotes an orthonormal basis of
$\ker B$. Consider the quadratic form,
\[
Q(\Phi)=\| D\Phi\| \quad \text{for} \quad \Phi\in
H^1_a(M_{1,\infty},S)
\]
where $D$ denotes the differential operator over $M_{1,\infty}$
whose self adjoint extension is $\Dd_{1,\infty}$. Then this
quadratic form is represented by a positive self adjoint operator
$H_a$ in the closure of $H^1_a(M_{1,\infty},S)$ in
$L^2(M_{1,\infty},S)$. Then $H_a$ has pure point spectrum near $0$
and $\ker H_a=\ker_{L^2} \Dd_{1,\infty}$ for any $a\ge 0$ by
Proposition 8.7 in \cite{Mu94}. Following the proof of Proposition
2.4 in \cite{KPW94}, we can prove that there exist positive
constants $c_1,c_2$ such that
\begin{equation}\label{orth}
| \langle \Phi_R, s\rangle | \le c_1 e^{-c_2R} \| s\|
\end{equation}
for $s\in\ker H_a$. Now let
$\tilde{\Phi}_R:=\Phi_R-\sum_{k=1}^{h_{1,\infty}}\langle\Phi_R,
s_k\rangle s_k$ where $\{s_k\}_{k=1}^{h_{1,\infty}}$ denotes an
orthonormal basis of $\ker H_a$ with $h_{1,\infty}:=\dim \ker
H_a$. Hence, $\tilde{\Phi}_R$ is orthogonal to $\ker H_a$, and by
\eqref{orth} there is a positive constant $C$ independent of $R$
such that $\|\tilde{\Phi}_R\| \ge \frac12 \|\Phi_R\| \ge C
> 0$ for sufficiently large $R$. Noting that $\tilde{\Phi}_R\in
\dom H_a$, and by the mini-max principle, we have
\begin{equation}\label{e:nu}
\langle H_a\tilde{\Phi}_R, \tilde{\Phi}_R\rangle \ge \nu^2 C^2
\end{equation}
where $\nu^2$ is the smallest nonzero eigenvalue of $H_a$. Now we
have
\begin{align*}
&\la(R)^2= \langle \Dd^2_R \Psi_R, \Psi_R\rangle = \int_{M_R} \|
\Dd_R \Psi_R(x)\|^2 \, dx \\ \ge&  \int_{M_{1,R}} \| \Dd_R
\Psi_R(x)\|^2 \, dx = \int_{M_{1,R}} \| \Dd_R \big(
h_R\Psi_R+(1-h_R)\Psi_R\big)(x)\|^2 \, dx\\
\ge& \int_{M_{1,\infty}} \| H_a\Phi_R(x)\| \, dx - \int_{M_{1,R}}
\| \Dd_R (1-h_R)\Psi_R(x)\|^2 \, dx \ \, .
\end{align*}
By \eqref{e:nu}, the first term has the lower bound $\nu^2C^2$
since $H_a\Phi_R=H_a\tilde{\Phi}_R$. For the second term, we have
\begin{align*}
 &\,\int_{M_{1,R}} \| \Dd_R (1-h_R)\Psi_R(x)\|^2 \, dx  \\=&\,
\int_{M_{1,R}} \|  (1-h_R)(x)\Dd_R \Psi_R(x)- G (\partial_u
h_R)(x)\Psi_R(x) \|^2 \, dx\\
\le &\, 2\int_{M_{1,R}} \|  \la(R)(1-h_R)(x)\Psi_R(x)\|^2+ \|G
(\partial_u h_R)(x)\Psi_R(x) \|^2 \, dx\ \, .
\end{align*}
By applying Lemma \ref{l:cut-est} with $v=u+R$ to each term of the
last line, we have
\[
 \int_{M_{1,R}}
\| \Dd_R (1-h_R)\Psi_R(x)\|^2 \, dx \le c_3 e^{-c_4 R}
\]
for positive constants $c_3,c_4$. Hence these inequalities imply
that $\la(R)^2\ge \frac12\nu^2C^2$ for sufficiently large $R$.
This completes the proof.

\end{proof}

Changing to the variable $v=u+R$, we regard that the cylindrical
part $N_R$ of $M_R$ is given by $[0,2R]_v\times Y$. In particular,
we have the new expression for $\Psi_R$ from \eqref{e:psi},
\begin{equation}\label{newexp}
\Psi_R = e^{-i\la v}\phi_1^1 + e^{i\la v}\phi_2^1 +\hat{\Psi}_R
\end{equation}
where $\phi_1^1=e^{i\la R}\psi_1$, $\phi_2^1=e^{-i\la R}\psi_2$.
Let $\Gg_{\pm}$ denote the $\pm i$ eigenspace of $G : \ker B \to
\ker B$. We need the following lemma,

\begin{lem}\label{l:inv} Let $\sigma$ be an involution over $\ker B$ such that
$G\sigma=-\sigma G$. Then for any element $\phi\in\Gg_\pm$, there
exists a unique $\psi\in\mathrm{Im}(\sigma+1)$ such that
\[
\phi=\psi\mp iG\psi.
\]
\end{lem}

\begin{proof}
For a given $\phi\in (\ker B)_+$, let
$\psi:=\frac12(1+\sigma)\phi$, which lies in
$\mathrm{Im}(\sigma+1)$ by definition. Then we have
\begin{align*}
\psi-iG\psi=&\frac12\big((1-i G)\phi + (\sigma-i G\sigma)\phi
\big)\\
=&\frac12\big((1-i G)\phi + (\sigma+i\sigma G)\phi \big) =
\frac12\cdot 2\phi=\phi\ \, .
\end{align*}
This completes the proof for the case of $(+)$ and the other case
of $(-)$ can be proved in the same way.
\end{proof}

By Proposition \ref{p:nts}, one of $\phi_1^1$ and $\phi_2^1$ in
\eqref{newexp} is nontrivial. First we assume that $\phi^1_1$ is
nontrivial. Now, since $L_1=\mathrm{Im}( C_1(0)+1)$ and $C_1(0)$
is an involution over $\ker B$, by Lemma \ref{l:inv} we can choose
$\psi\in L_1$ such that $\phi_1^1=\psi-iG\psi$. Then the
generalized eigensection $E(\psi,\la)$ over $M_{1,\infty}$
associated to $\psi$ has the following expression
$$
E(\psi,\lambda)=e^{-i\lambda v}(\psi-iG\psi)+ e^{i\la
v}C_{1}(\la)(\psi-iG\psi) + \theta(\psi,\la)$$ over $[  0,
\infty)_v\times Y$ . Following \cite{Mu94}, we introduce
$$
F=\Psi_R|_{M_{1,R}} - E(\psi,\lambda)|_{M_{1,R}}\ \, . $$ Green's
formula gives
$$
0 = \langle \Dd_{1,R} F, F\rangle_{M_{1,R}} - \langle F, \Dd_{1,R}
F \rangle_{M_{1,R}} = \int_{\partial{(M_{1,R})}} \langle G
F|_{\partial{(M_{1,R})}}, F|_{\partial{(M_{1,R})}} \rangle\, dy.$$
On the other hand, Lemma \ref{l:cut-est} shows that
$$
\int_{\partial{(M_{1,R})}} \langle G F|_{\partial{(M_{1,R})}},
F|_{\partial{(M_{1,R})}} \rangle\, dy= -i\|\ C_1(\lambda)\phi_1^1
-\phi_2^1 \|^2 + O(e^{-c_3 R})
$$
for some positive constant $c_3$. This produces the estimate
\begin{equation}\label{cc1}
\|\, C_1(\la)\phi_1^1 - \phi_2^1\, \|  \le e^{-cR} \ \,
\end{equation}
for a positive constant $c$. Therefore, for $R\gg 0$, if
$\phi^1_1$ is nontrivial, then $\phi^1_2$ is also nontrivial. In
the same way, one can show its inverse. Hence we can conclude that
both $\phi^1_1,\phi^1_2$ in \eqref{newexp} are nontrivial for
$R\gg 0$.

Now we want to get the corresponding estimate involving the
scattering matrix $C_2(\la)$. For this, we change the variable by
$v=u-R$ and regard the cylindrical part as $[-2R,0]_v\times Y$.
Then we have the corresponding expression for $\Psi_R$,
\[
\Psi_R = e^{-i\la v}\phi_1^2 + e^{i\la v}\phi_2^2 +\hat{\Psi}_R
\]
where $\phi_1^2=e^{-i\la R}\psi_1$, $\phi_2^2=e^{i\la R}\psi_2$.
For the given $\phi^2_2\in \Gg_-$, using Lemma \ref{l:inv}, we
choose $\psi\in L_2=\mathrm{Im}(C_2(0)+1)$ such that
$\phi^2_2=\psi+iG\psi$. The generalized eigensection $E(\psi,\la)$
over $M_{2,\infty}$ attached to the couple $(\psi,\la)$ has the
following expression
$$
E(\psi,\lambda)=e^{i\lambda v}(\psi+iG\psi)+ e^{-i\la
v}C_{2}(\la)(\psi+iG\psi) + \theta(\psi,\la)$$ over $(-
\infty,0]_v\times Y$. As above, comparing $\Psi_R$ and
$E(\psi,\la)$, we obtain
\begin{equation}\label{cc2}
\|\, C_2(\la)\phi_2^2 - \phi_1^2\, \|  \le e^{-cR} \ \,
\end{equation}
for a positive constant $c$.  By definition, we have
\begin{equation}\label{e:phi}
\phi^1_1=e^{2i\la R}\phi^2_1\ \ , \quad  \phi^1_2=e^{-2i\la
R}\phi^2_2\ \ .
\end{equation}
Now, combining \eqref{cc1}, \eqref{cc2} and \eqref{e:phi}, we get
\begin{equation}\label{c12}
\|\, e^{4i\la R} C_1(\la)\circ C_2(\la)\phi_2^1 - \phi_2^1\, \|
\le e^{-cR} \ \ .
\end{equation}

We define the operator $C_{12}(\lambda)$ by
$$
C_{12}(\lambda) := C_1(\lambda)\circ C_2(\lambda)|_{\Gg_- } \ :\
\Gg_{-} \to \Gg_{-} \ \ .
$$
The operator $C_{12}(\lambda)$ is a unitary operator and is an
analytic function of $\lambda$ for $\lambda\in (-\delta,\delta)$
for a small $\delta>0$ since the unitary operators $C_1(\lambda)$,
$C_2(\lambda)$ are analytic functions of $\lambda$ for $\lambda
\in (-\delta, \delta)$. Furthermore, there exist real analytic
functions $\alpha_j(\lambda)$ for $1\le j\le \frac{h_Y}2$ of
$\lambda\in(-\delta,\delta)$ such that $\exp(i\alpha_j(\lambda))$
are the corresponding eigenvalues of $C_{12}(\lambda)$ and
$\alpha_j(\la)$ has the following expansion at $\la=0$,
\begin{equation}\label{e:expansion}
\alpha_j(\la)=\alpha_{j0}+\alpha_{j1}\la+\alpha_{j2}\la^2+\alpha_{j3}\la^3+\ldots\
\, .
\end{equation}
We now introduce
\begin{equation}\label{OmegaR}
\Omega(R):=\{\, \rho\in\mathbb{R}-\{0\} \ | \ \det(e^{4i\rho
R}C_{12}(\rho)-\Id)=0 \ ,\ |\rho| \le R^{-\kappa}\, \} \ \, .
\end{equation}
The following theorem is a main result of this section,

\begin{thm}\label{t:small eigen1} Assume that all the
\emph{e-values} of $\Dd_R$ are zero eigenvalues. Let $\la_1(R)\le
\la_2(R)\le \ldots \le \la_{p(R)}(R)$ be the nonzero eigenvalues,
counted to multiplicity, of $\Dd_R$ satisfying $|\la_k(R)|\le
R^{-\kappa}$, and let $\rho_1(R)\le \rho_2(R)\le \ldots \le
\rho_{m(R)}(R)$ be the nonzero element , counted to multiplicity,
of $\Omega(R)$. Then there exist $R_0$ and $c>0$, independent of
$R$,  such that for $R\ge R_0$, $p(R)=m(R)$ and
\[
|\la_k(R) -\rho_k(R)| \le e^{-cR} \quad \text{for}\quad
k=1,\ldots, p(R)\ \, .
\]
\end{thm}

\begin{proof}
The proof of this theorem consists of two steps. {\flushleft Step
I: } Let $\la=\la(R)$ be a given \emph{s-value} with the
multiplicity $m(\la)$. By Proposition \ref{p:nts}, we get $m(\la)$
linearly independent vectors $\phi_1, \ldots, \phi_{m(\la)}$ in
$\Gg_-$, which satisfies \eqref{c12}. Since $C_{12}(\la)$ is
unitary, the eigenvalues of $e^{4i\la R }C_{12}(\la)-\Id$ have the
form $e^{i\theta}-1$ for $\theta\in\mathbb{R}$. Let $0\le \z$ be
the smallest eigenvalue of
 $(e^{4i\la R }C_{12}(\la)-\Id)(e^{4i\la R }C_{12}(\la)-\Id)^*$;
then
\[
\z=\min_{\phi\in \Gg_-} \frac{\|(e^{4i\la
R}C_{12}(\la)-\Id)\phi\|^2}{\|\phi\|^2} \ \, .
\]
Combined with \eqref{c12}, this implies that $\z\le e^{-cR}$.
Hence $e^{4i\la R}C_{12}(\la)$ has an eigenvalue $e^{i\theta}$
satisfying $|1-\cos\theta| \le e^{-cR}$, and there exists
$k\in\mathbb{Z}$ such that $|2\pi k -\theta| \le e^{-cR}$.
Therefore, by definition of $\alpha_j(\la)$, the following holds
\begin{equation}\label{equ}
|\, 4\lambda R + \alpha_j(\lambda) -2\pi k \, |\le e^{-cR}
\end{equation}
for pairwise distinct branches $\alpha_1, \ldots,
\alpha_{m(\la)}$. Now, let us fix $\delta_1$ with
$0<\delta_1<\delta$ and let
\begin{equation*}
m_j=\max_{\la\in (-\delta_1,\delta_1)} |\alpha_j'(\la)| \ \, .
\end{equation*}
Then the function $f(\la)=4\la R + \alpha_j(\la)$ is strictly
increasing for $|\la|<\delta_1$ and $R\ge m_j$.  Choose $R_1$ such
that $R_1 \ge \max(m_j,\delta_1^{-\frac1\kappa})$ for any
$j=1,\dots,\frac{h_Y}2$. For $R\ge R_1$ and $k\in\mathbb{Z}$,
there exists at most one solution $\rho_{j,k}$ of
\begin{equation}\label{e:sol}
4 \la R +\alpha_j(\la) = 2\pi k \ , \qquad |\la| \le R^{-\kappa} \
\, .
\end{equation}
Let $k_{j,\max}$ be the maximal $k$ for which \eqref{e:sol} has a
solution; then by \eqref{e:sol},
\begin{equation}\label{e:kmax}
|k_{j,\max}| \le  \frac{2R^{1-\kappa}}{\pi} +C \le R^{1-\kappa} \
\,.
\end{equation}
Then, for $R\ge R_1$, any element $\rho$ in $\Omega(R)$ is given
by $\rho=\rho_{j,k}$ for some $1\le j \le \frac{h_Y}2$, and
$|k|\le k_{j,\max}$. Therefore, if $R\ge R_1$, for a given $\la$
satisfying \eqref{equ} with $|\la| \le R^{-\kappa}$, there is a
unique solution $\rho_{j,k}$ of \eqref{e:sol} such that
\begin{equation}\label{erelation}
|\la-\rho_{j,k}| \le e^{-cR}\  \, .
\end{equation}
In conclusion, if $R\ge R_1$, for a given \emph{s-value}
$\la=\la(R)$ of $\Dd_R$ with the multiplicity $m(\la)$ satisfying
$|\la|\le R^{-\kappa}$, there exist $m(\la)$-number of elements
$\rho_{j,k}$'s in $\Omega(R)$ with the relation \eqref{erelation},
in particular, $p(R)\le m(R)$.

{\flushleft Step II:} To complete the proof, we need to prove that
$m(R)\le p(R)$. For $k$ with $1\le k \le m(R)$, we choose
$\psi_k\in \Gg_-$ with the following properties,
\begin{enumerate}
\item  $e^{4i\rho_k R}C_{12}(\rho_k)\psi_k =\psi_k \ , \quad
|\rho_{k}|\le R^{-\kappa}\ \, , $ \item When
$\rho_k=\rho_{k+1}=\ldots =\rho_{k+\ell}$,
$\psi_k,\psi_{k+1}\ldots, \psi_{k+\ell}$ form an orthonormal
system of vectors of $\Gg_-$\ .
\end{enumerate}
For a given pair $(\psi_k, \rho_k)$ for some $k$, we put
\begin{equation}\label{connecting}
\phi^1_k=e^{-i\rho_k R}C_1(-\rho_k)\psi_k\, , \quad
\phi^2_k=e^{i\rho_k R}\psi_k\ \, .
\end{equation}
Now we consider the generalized eigensection
$E({\phi}^1_k,\rho_k)$ over $M_{1,\infty}$ and
$E({\phi}^2_k,\rho_k)$ over $M_{2,\infty}$, which have the
following forms,
\begin{align*}
E(\phi^1_k,\rho_k)\ =& \ e^{-i\rho_k v}\phi^1_k+e^{i\rho_k
v}C_1(\rho_k)\phi^1_k+ O(e^{-cv}) \quad \text{over}\quad
[0,\infty)_v\times Y \subset
M_{1,\infty} \ \, ,\\
E(\phi^2_k,\rho_k)\ =&\  e^{i\rho_k v}\phi^2_k+e^{-i\rho_k
v}C_2(\rho_k)\phi^2_k+ O(e^{-cv}) \ \ \text{over}\ \
(-\infty,0]_v\times Y \subset M_{2,\infty}\ \, .
\end{align*}
(Here we use abuse notations for simplicity since the correct
notation for $E(\phi^i_k,\rho_k)$ is $E(\tilde{\phi}^i_k,\rho_k)$
with $\phi^i_k=\tilde{\phi}^i_k+ (-1)^i
\sqrt{-1}G\tilde{\phi}^i_k$ by Lemma \ref{l:inv}.) Restricting
$E(\phi^i_k,\rho_k)$ to $M_{i,R}$, we obtain sections over
$M_{i,R}$. Let $f_{1,R}$ be the restriction to $M_{1,R}$ of the
smooth function $h_R$ over $M_{1,\infty}$ defined in the proof of
Proposition \ref{p:nts} and $f_{2,R}$ be a smooth function over
$M_{2,R}$ defined in a similar way.  These functions have the
obvious extension over $M_R$.  Denoting by $E_0(\phi^i_k,\rho_k)$
the zero eigenmode of $E(\phi^i_k,\rho_k)$ and using
\eqref{connecting} and $e^{4i\rho_k R}C_{12}(\rho_k)\psi_k =\psi_k
$, we have
\begin{multline}\label{const}
E_0(\phi^1_k,\rho_k)=e^{-i\rho_k u}
e^{-2i\rho_kR}C_1(-\rho_k)\psi_k+e^{i\rho_k u}\psi_k \\
=e^{i\rho_k u}\psi_k +e^{-i\rho_k u}
e^{2i\rho_kR}C_2(\rho_k)\psi_k= E_0(\phi^2_k,\rho_k)\ \, .
\end{multline}
Hence we can see that $E_0(\phi^1_k,\rho_k)$ and
$E_0(\phi^2_k,\rho_k)$ define a smooth section over $N_R$, which
we denote by $E_0(\psi_k,\rho_k)$. Let us define
\begin{multline}\label{tildepsi}
\tilde{\Psi}_k:= f_{1,R}(E(\phi^1_k,\rho_k)
-\chi_{[-R,0]_u}E_0(\phi^1_k,\rho_k)) \\ +f_{2,R}
(E(\phi^2_k,\rho_k)-\chi_{[0,R]_u}E_0(\phi^2_k,\rho_k))
+\chi_{[-R,R]_u} E(\psi_k,\rho_k)
\end{multline}
where $\chi_{[a,b]_u}$ is the characteristic function of the
$u$-variable over $[a,b]_u\times Y\subset N_R$. By \eqref{const},
$\tilde{\Psi}_k$ is a smooth section over $M_R$. Put
$\Psi_k:=\tilde{\Psi}_k/ \|\tilde{\Psi}_k\|$ and
\[
\hat{\Psi}_k=\Psi_k -\pi_R \Psi_k, \quad\text{for} \quad
k=1,\ldots, m(R) \ \, ,
\]
where $\pi_R$ denote the orthogonal projection of $L^2(M_R,S)$
onto $\ker \Dd_R$. Let us recall that $\ker \Dd_R$ equals the
space spanned by eigensections of \emph{e-values} by our
assumption, so that the dimension of this space is constant with
respect to $R$. Combining this fact and Lemma \ref{orthogonal}, we
have
\[
| \langle \hat{\Psi}_k, \hat{\Psi}_\ell \rangle -\delta_{k\ell}|
\le e^{-cR} \ \,  \quad\text{for}\quad k,\ell=1,\ldots, m(R) \ \,
.
\] From this and \eqref{e:kmax}, it follows that
$\{\hat{\Psi}_k\}_{k=1}^{m(R)}$ are linearly independent for $R\gg
0$. Now let $0 <\tilde{\la}_1 \le \tilde{\la}_2\le \ldots\le
\tilde{\la}_{p(R)}$ denote the nonzero eigenvalues, counted with
multiplicity, of $\Dd_R^2$, which are $\le R^{-2\kappa}$. Let
$k_1,\ldots, k_{m(R)}$ be a permutation of $\{1,\ldots, m(R)\}$
such that $0< \rho^2_{k_1}\le \ldots \le \rho^2_{k_{m(R)}}$. By
the mini-max principle, we have
\[
\tilde{\la}_{\ell} = \min_W \max_{\phi\in W} \frac{\| \Dd_R \phi
\|^2}{\|\phi\|^2} \ \,
\]
where $W$ runs over all $\ell$-dimensional subspaces of
$L^2(M_R,S)$ which are orthogonal to $\ker(\Dd_R)$. Let $W_{\ell}$
be the subspace of $L^2(M_R,S)$ spanned by $\hat{\Psi}_{k_1},
\ldots, \hat{\Psi}_{k_{\ell}}$. Then, by Lemma \ref{orthogonal},
we have
\[
\tilde{\la}_{\ell} \le \max_{\phi\in W_\ell} \frac{\| \Dd_R \phi
\|^2}{\|\phi\|^2} \le \rho^2_{k_{\ell}} (1+ C e^{-cR})
\]
for some constants $C,c>0$. Hence, there exists $R_2$ such that
$m(R)\le p(R)$ for $R\ge R_2$. Putting $R_0=\max(R_1,R_2)$, this
completes the proof of Theorem \ref{t:small eigen1}.
\end{proof}

\begin{lem}\label{orthogonal} Assume that all the
\emph{e-values} of $\Dd_R$ are zero eigenvalues. Then there exist
$c_1,c_2>0$ such that
\begin{align*}
|\langle \Psi_k, \Psi_\ell \rangle | \le c_1 e^{-c_2 R}& \quad
\text{for}\quad k\neq \ell, \ k,\ell=1,\ldots, m(R) \ \, ,\\
|\langle \Psi_k, \Psi \rangle | \le c_1 e^{-c_2 R}& \quad
\text{for}\quad  k =1,\ldots, m(R) , \ \text{and} \ \Psi\in \ker
\Dd_R \ \text{with} \ ||\Psi||=1 \ \, .
\end{align*}
\end{lem}

\begin{proof}
For a couple $(\psi_k,\rho_k)$ and $\phi^i_k$ satisfying
\eqref{connecting}, we put
\begin{align*}
E^\perp_k=& E(\phi^1_k,\rho_k)|_{M_{1,R}}
-\chi_{[-R,0]_u}E_0(\phi^1_k,\rho_k) + E(\phi^2_k,\rho)|_{M_{2,R}}
-\chi_{[0,R]_u}E_0(\phi^2_k,\rho_k)\ \, ,\\
E_{k,0}=& E(\psi_k,\rho_k) = \chi_{[-R,0]_u}E_0(\phi^1_k,\rho_k)
+\chi_{[0,R]_u}E_0(\phi^2_k,\rho_k) \ \, .
\end{align*}
Putting $f_R=f_{1,R}+f_{2,R}$, it is easy to see that
$\tilde{\Psi}_k$ defined in \eqref{tildepsi} has the form $f_R
E_k^\perp + \chi_{[-R,R]_u}E_{k,0}$. Now we have
\begin{align}
&\ \langle \tilde{\Psi}_k, \tilde{\Psi}_{\ell} \rangle = \langle f
E^\perp_k +\chi
E_{k,0}, f E^\perp_{\ell} +\chi  E_{\ell,0} \rangle \notag \\
=&\ \langle f E^\perp_k , f E^\perp_{\ell} \rangle +\langle \chi
E_{k,0}, \chi  E_{\ell,0} \rangle \notag \\
=&\ \langle  E^\perp_k -(1-f)E^\perp_k,
E^\perp_{\ell}-(1-f)E^\perp_{\ell} \rangle +\langle \chi
E_{k,0}, \chi  E_{\ell,0} \rangle \notag \\
=&\ \langle E^\perp_k ,E^\perp_{\ell} \rangle -\langle E^\perp_k,
(1-f)E^\perp_\ell \rangle - \langle (1-f) E^\perp_k, E^\perp_\ell
\rangle \notag \\ &+ \langle (1-f) E^\perp_k, (1-f)E^\perp_\ell
\rangle +\langle \chi E_{k,0}, \chi E_{\ell,0}
\rangle \notag \\
=&\ \langle E_k ,E_{\ell} \rangle -\langle E^\perp_k,
(1-f)E^\perp_\ell \rangle - \langle (1-f) E^\perp_k, E^\perp_\ell
\rangle \label{threeterm}\\ &+ \langle (1-f) E^\perp_k,
(1-f)E^\perp_\ell \rangle \notag
\end{align}
where $f=f_R, \chi=\chi_{[-R,R]_u}$. Since $\supp(1-f_R)\subset
[-\frac R2, \frac R2]_u\times Y$, where $E^\perp_k, E^\perp_\ell$
are $O(e^{-cR})$, the last three terms in \eqref{threeterm} are
$O(e^{-cR})$. Now we consider the first term in \eqref{threeterm},
which can be written as
\begin{equation}\label{E}
\langle E_k ,E_{\ell} \rangle = \langle E(\phi^1_k,\rho_k),
E(\phi^1_\ell,\rho_\ell) \rangle_{M_{1,R}} + \langle
E(\phi^2_k,\rho_k), E(\phi^2_\ell,\rho_\ell) \rangle_{M_{2,R}}\ \,
.
\end{equation}
When $\rho_k\neq \rho_\ell$, as in the proof of Proposition
\ref{Maass-Selberg}, we apply Green formula to each term on the
right side of \eqref{E}, then these equal
\begin{multline*}
 (\rho_k-\rho_\ell)^{-1}\langle G E(\phi^1_k,\rho_k)|_{\partial(M_{1,R})},
E(\phi^1_\ell,\rho_\ell)|_{\partial(M_{1,R})}
\rangle_{\partial(M_{1,R})} \\- (\rho_k-\rho_\ell)^{-1}\langle G
E(\phi^2_k,\rho_k)|_{\partial(M_{2,R})},
E(\phi^2_\ell,\rho_\ell)|_{\partial(M_{2,R})}
\rangle_{\partial(M_{2,R})} \ \, .
\end{multline*}
Now using \eqref{const}, the restrictions of constant terms over
$\partial(M_{i,R})$ cancel each other out and the remaining terms
are $O(e^{-cR})$. Hence, in this case, the left side of \eqref{E}
is $O(e^{-cR})$, so all the terms in \eqref{threeterm} are
$O(e^{-cR})$. When $\rho_k=\rho_l$, note that $\langle \phi^i_k,
\phi^i_\ell \rangle =0 $ for $i=1,2$, so applying Proposition
\ref{Maass-Selberg}, we can see that all the terms are
$O(e^{-cR})$ except the following terms,
\begin{equation}\label{non}
\langle C_1(-\rho_k) C_1'(\rho_k)\phi^1_k, \phi^1_\ell \rangle +
 \langle C_2(-\rho_k) C_2'(\rho_k)\phi^2_k, \phi^2_\ell \rangle \
 \, .
\end{equation}
Using the conditions in \eqref{connecting} for
$\phi^i_k,\phi^i_\ell$ and  the relation
\begin{equation}\label{relation}
e^{4i\rho_k R} C_2(\rho_k)\psi_k = C_1(-\rho_k)\psi_k \ , \quad
e^{4i\rho_k R} C_2(\rho_k)\psi_\ell = C_1(-\rho_k)\psi_\ell\ \, ,
\end{equation}
one can show that the terms in \eqref{non} equal
\begin{equation}\label{EE}
\langle e^{4i\rho_k R } C_1'(\rho_k)C_2(\rho_k)\psi_k, \psi_\ell
\rangle + \langle e^{4i\rho_k R} C_1(\rho_k)C_2'(\rho_k)\psi_k,
\psi_\ell \rangle\ \, .
\end{equation}
Now we choose a family of sections $\psi_k(t)$ with
$\psi_k(0)=\psi_k $ for $t\in(-\epsilon,\epsilon)$ such that
$$
a(t)C_1(\rho_k+t) C_2(\rho_k+t)\psi_k(t) = \psi_k(t)
$$
where $a(0)=e^{4i\rho_k R }$. Taking the derivative of this at
$t=0$, we obtain
\begin{multline*}
e^{4i\rho_k R } C_1'(\rho_k)C_2(\rho_k)\psi_k +  e^{4i\rho_k R }
C_1(\rho_k)C_2'(\rho_k) \psi_k\\
= -{a}'(0) C_1(\rho_k)C_2(\rho_k)\psi_k -  e^{4i\rho_k R }
C_1(\rho_k)C_2(\rho_k){\psi}'_k(0) +{\psi}'_k(0) \ \, .
\end{multline*}
Using this, \eqref{relation} and $\langle \psi_k,\psi_\ell
\rangle=0$, we can see that \eqref{EE} equals
\begin{multline*}
\langle  -  e^{4i\rho_k R } C_1(\rho_k)C_2(\rho_k){\psi}'_k(0)
+{\psi}'_k(0) , \psi_\ell \rangle\\
= \langle {\psi}'_k(0), \psi_{\ell} \rangle - \langle
{\psi}'_k(0), e^{-4i\rho_k R } C_2(-\rho_k)C_1(-\rho_k)\psi_\ell
\rangle =0 \ \, .
\end{multline*}
Hence, in the case of $\rho_k=\rho_l$, all the terms in
\eqref{threeterm} are $O(e^{-cR})$. This completes the proof of
the first claim recalling
$\Psi_k=\tilde{\Psi}_k/||\tilde{\Psi}_k||$.

For the second claim, let us recall that the eigenspaces of the
\emph{e-values} are spanned by the sections defined by gluing (as
in \eqref{tildepsi}) the elements in $\ker_{L^2} \Dd_{i,\infty}$
for $i=1,2$ or the extended $L^2$-solutions of $\Dd_{i,\infty}$
whose limiting values lying in $L_1\cap L_2$. By our assumption,
this space is the same as $\ker \Dd_R$. For a section $\Psi$ given
by gluing elements in $\ker_{L^2}\Dd_{i,\infty}$, the claim
follows easily by applying Green's formula as above. For a section
$\Psi$ given by gluing the extended $L^2$-solutions of
$\Dd_{i,\infty}$ whose limiting values lying in $L_1\cap L_2$, we
use Theorem \ref{t:eev}, which implies that such a $\Psi$ is
actually given by \eqref{tildepsi} for the couple
$(\psi_k,\rho_k)$ with $\rho_k=0$. Hence, the claim for this case
can be proved as in the previous case of $\rho_k\neq\rho_\ell$.
This completes the proof of the second claim.

\end{proof}

In general, the map $C_{12}:=C_{12}(0): \Gg_- \to \Gg_-$ does not
equal the identity map, but it is not difficult to see that
$$
C_1(0)\circ C_2(0) \phi=\phi \quad \text{if and only if} \quad
\phi\in (L_1\cap L_2) \oplus (GL_1\cap GL_2)\ \, .
$$
Putting
$$
I_+= 1+iG:\ker B \to \Gg_- \ \ ,
$$
we can see that $I_+(L_1\cap L_2)$ and $I_+(GL_1\cap GL_2)$ are
the same subspace in $\Gg_-$.

\begin{prop}\label{p:L12}
The map $C_{12}$ equals the identity map when restricted to the
subspace $I_+(L_1\cap L_2)$ and the multiplicity of the eigenvalue
$(+1)$ of the operator $C_{12}$ is  $\dim(L_1\cap L_2)=\dim
(I_+(L_1\cap L_2))$.
\end{prop}

\begin{proof}
Using the following diagram
$$
\begin{CD}
L_1\cap L_2 @>I_+ >> \Gg_-\\
@V
C_1(0)\circ C_2(0)VV @VVC_1(0)\circ C_2(0)V \\
L_1\cap L_2 @>I_+ >> \ \ \Gg_- \ \, ,
\end{CD}
$$
we can easily see that the first claim holds. To complete the
proof, it is sufficient to show that if $C_{1}(0)\circ
C_2(0)\phi=\phi$, then $\phi \in (L_1\cap L_2) \oplus (GL_1\cap
GL_2)$. For this, choose $\phi_+\in L_1$, then $C_{1}(0)\circ
C_2(0) \phi_+=\phi_+$ implies $\phi_+=C_1(0)\phi_+=C_2(0)\phi_+$
since $C_1(0)^2=\Id$. Hence, this means that $\phi_+\in L_2$, so
$\phi_+\in L_1\cap L_2$. Repeating the same argument, if
$\phi_-\in GL_1$ and $C_{1}(0)\circ C_2(0)\phi_-=\phi_-$, then
$\phi_-\in GL_1\cap GL_2$. Since $\ker B=L_1\oplus GL_1$, this
completes the proof.

\end{proof}

Now let us consider the eigenvalues $\la(R)$ , which correspond to
$\alpha_j(0) = 0$ and $k=0$ in the following equality equivalent
to \eqref{equ},
\begin{equation*}
4\lambda R + \alpha_j(\lambda) =2\pi k + O(e^{-cR})\ \ .
\end{equation*}
It is easy to see that such eigenvalues must be \emph{e-values}.
Hence, by Lemma \ref{p:L12} this provides another proof of the
following result, originally shown in \cite{CLM1}.

\begin{thm}\label{t:eev}
The space of eigensections corresponding to \emph{e-values}, which
are not determined by $\ker_{L^2}(\Dd_{i,\infty})$ for $i=1,2$, is
given by the space $L_1\cap L_2$.
\end{thm}

Now let us consider the following Dirac type operator
\begin{equation}\label{hodge}
\Dd_R=d_{\rho}+d_{\rho}^* \ : \ \oplus_{i=0}^{n} \Omega^i(M_R,
V_{\rho}) \ \rightarrow \ \oplus_{i=0}^{n} \Omega^i(M_R, V_{\rho})
\end{equation}
where $V_{\rho}$ denotes the flat vector bundle defined by a
unitary representation $\rho$ of $\pi_1(M_R)$. The dimension of
$\ker \Dd_R$ is constant with respect to $R$ since $\ker \Dd_R$ is
the space of the twisted harmonic forms over $M_R$ and this space
is always isomorphic to de Rham cohomology $H^*(M_R,V_\rho)$ by
the Hodge theorem. Moreover, one can show that all the
\emph{e-values} of the operator $\Dd_R$ in \eqref{hodge} are the
zero eigenvalues using the argument in Section 4 of \cite{Ha98}.

\begin{prop} \label{p:zero} For the operator $\Dd_R$ in \eqref{hodge},
all the \emph{e-values} of $\Dd_R$ are the zero eigenvalues.
\end{prop}

\begin{proof}
First let us observe that $\Dd_{i,\infty}$ is self adjoint, so
$\ker_{L^2}\Dd_{i,\infty}=\ker_{L^2}\Dd_{i,\infty}^2$ and $L_i$ is
also the limiting value of extended $L^2$-solutions of
$\Dd^2_{i,\infty}$. Let $\Del_{i,\infty}^q$ be the restriction of
$\Dd_{i,\infty}^2$ to $\Omega^q(M_R, V_{\rho})$ and
\begin{equation}\label{defh}
 h_M^q:=
\dim \ker_{L^2} \Del_{1,\infty}^q + \dim \ker_{L^2}
\Del_{2,\infty}^q + \dim L_1^q\cap L_2^q \end{equation} where
$L_i^q$ is the limiting values of the extended $L^2$-solutions of
$\Delta^q_{i,\infty}$. Then, it is sufficient to show that
$\beta^q:=\dim (\ker \Dd^2_R\cap \Omega^q(M_R, V_{\rho})) \ge
h_M^q$ since $\beta_q\le h_M^q$ by definition. For this, we use
the following Mayer-Vietoris sequence
\begin{equation}\label{MV}
\ldots \to H^{q-1}(Y)\to H^q(M_R) \to H^q_a(M_{1,R})\oplus
H^q_a(M_{2,R}) \to H^q(Y) \to \ldots \ \,
\end{equation}
where $H^q_a(M_{i,R})$ denotes the absolute cohomology. (Here, for
simplicity, the bundle $V_\rho$ is dropped in the notation.) The
space $\oplus_{q=0}^{n} H^q_a(M_{i,R})$ can be identified with the
kernel of the operator $\Dd_{i,R}$ with the absolute boundary
condition. In more detail, the operator $\Dd_R=d_\rho+d^*_\rho$
has the following form over $N_R$,
\begin{equation}\label{matrix}
\Dd_R=d_\rho+d_\rho^*= \begin{pmatrix} 0 & -1 \\ 1 & 0
\end{pmatrix} \Big(\partial_u +\begin{pmatrix} 0 & d_Y+d_Y^* \\
d_Y+d_Y^* & 0
\end{pmatrix} \Big)
\end{equation}
with respect to
\begin{equation}\label{decomp}
\Omega^*(N_R) \cong (\Omega^*(Y)\oplus \Omega^*(Y))\otimes
C^{\infty}([-R,R]_u).
\end{equation}
Here $d_Y$, $d_Y^*$ denote the restricted operator to $Y$ of
$d_\rho$ and its adjoint respectively. The operator $\Dd_{i,R}$
has the same form near the boundary and with respect to
\eqref{decomp}. A section $\Psi$ in $\Omega^*(M_{i,R})$ over the
cylinder near the boundary $Y$ has the following form,
\[
\Psi = \Psi_0 +\Psi_1\wedge du \ \,
\]
where $\Psi_i$ has no factor $du$. Then the absolute boundary
condition for $\Dd_{i,R}$ is given by $\Psi_1=0$. Similarly the
relative boundary condition for $\Dd_{i,R}$ is given by
$\Psi_0=0$. We denote by $\Dd_{i,R}^a$, $\Dd_{i,R}^r$ the
resulting operators. Now let us recall that the Cauchy data spaces
$\Hh(\Dd_{i,R})$ of $\Dd_{i,R}$ are Lagrangian subspaces in
$\Omega^*(Y)\oplus \Omega^*(Y)$ with respect to the symplectic
form $\langle G \, , \, \rangle$ where $G=\begin{pmatrix} 0 & -1 \\
1 & 0
\end{pmatrix}$ and $\langle \, , \, \rangle$ are the inner product
over $\Omega^*(Y)\oplus \Omega^*(Y)$. Then, this implies that $
\Hh_0(\Dd_{i,R}):=\Hh(\Dd_{i,R})\cap (H^*(Y)\oplus H^*(Y))$ are
also Lagrangian subspaces in $H^*(Y)\oplus H^*(Y)$. Moreover, the
space $\Hh_0(\Dd_{i,R})$ has the following decomposition,
\[
\Hh_0(\Dd_{i,R})=A_i\oplus R_i
\]
where $A_i,R_i$ the spaces spanned by the boundary values of $\ker
\Dd_{i,R}^a, \ker \Dd_{i,R}^r$ in $H^*(Y)\oplus H^*(Y)$.
Decomposing  $A_i=\oplus_{q=0}^{n-1} A^q_i$, $R_i=
\oplus_{q=0}^{n-1} R^q_i$  where $A^q_i,R^q_i$ are spaces of
$q$-form parts, the Lagrangian subspace property of
$\Hh_0(\Dd_{i,R})$ in $H^*(Y)\oplus H^*(Y)$ implies
\begin{equation}\label{e:decompH}
H^q(Y) \cong A^q_i \oplus R^{q}_i\ \, .
\end{equation}
By the exactness of \eqref{MV}, we also have
\begin{align*}
H^q(M_R) &\cong \mathrm{Im}(H^{q-1}(Y)\to H^q(M_R))\oplus
\mathrm{Im}(H^q(M_R) \to H^q_a(M_{1,R})\oplus H^q_a(M_{2,R})) \\
&\cong (\mathrm{Im}\, k^{q-1})^\perp \oplus \ker k^q \ \,
\end{align*}
where $k^q$ is the boundary map from $H^q_a(M_{1,R})\oplus
H^q_a(M_{2,R})$ to $H^q(Y)$. Now we summarize the consequences of
the previous considerations. First, by \eqref{e:decompH}, we have
\[
(\mathrm{Im}\, k^{q-1})^\perp = (A_1^{q-1}+A_2^{q-1})^\perp =
R^{q-1}_1\cap R^{q-1}_2 \ \, .
\]
Second, we note that $\ker k^q$ contains the harmonic sections
whose boundary values are lying in $A^q_1\cap A^q_2$ and the
harmonic sections that can be extended as $L^2$-solutions of
$\Del^q_{i,\infty}$. Hence,
\[
\dim \ker k^q \ge \dim(A^q_1\cap A^q_2)+\dim \ker_{L^2}
\Del_{1,\infty}^q + \dim \ker_{L^2} \Del_{2,\infty}^q\ \, .
\]
By these facts and the following equality
\[
\dim(L^q_1\cap L^q_2) = \dim(A^q_1\cap A^q_2)+\dim(R^{q-1}_1\cap
R^{q-1}_2), \] we can conclude that $\beta^q \ge h_M^q$ recalling
\eqref{defh}. This completes the proof.

\end{proof}

Let $\Psi_R$ be a normalized eigensection of $(\Dd_{1,R})_{P_1}$ ,
which corresponds to the \emph{s-value} $\la = \la(R)$ with
$|\lambda(R)| \le R^{-\kappa}$ where $\kappa$ is the fixed
constant such that $0<\kappa<1$ . Then we have
\begin{equation}\label{e:boundary}
\Dd_{1,R}\Psi_R=\lambda\Psi_R  \ \ , \ \ \|\Psi_R\| =1 \ \, , \ \
P_1(\Psi_R|_{\{v=R\}\times Y})=0 \ \ .
\end{equation}
The section $\Psi_R$ can be represented in the following way on
$[0,R]_v\times Y \subset M_{1,R}$
\begin{equation}\label{express}
\Psi_R=e^{-i\lambda v}\psi_1 + e^{i\lambda v}\psi_2
+{\hat{\Psi}_R}
\end{equation}
where $\psi_1\in \Gg_+ , \psi_2 \in \Gg_-$ and ${\hat{\Psi}_R}$ is
a smooth $L^2$-section orthogonal to $\ker B$. The next result
corresponds to Proposition \ref{p:nts}, which can be proved in the
same way as Proposition 8.14 of \cite{Mu94},

\begin{prop}\label{p:nts1}
The zero-eigenmode $e^{-i\lambda v}\psi_1+ e^{i\lambda v}\psi_2$
of the eigensection $\Psi_R$ of \emph{s-value} $\la(R)$ of
$(\Dd_{1,R})_{P_1}$ is a non-trivial.
\end{prop}

Now we define
\begin{align*}
I_\pm&= 1 \pm iG : \ker B \to \Gg_\mp \ \ ,\\
I_{\sigma_1} &= I_-|_{\ker(\sigma_1+1)} : \ker (\sigma_1+1) \to \Gg_+ \ \ ,\\
P_{\sigma_1}&=\frac 12(1- \sigma_1) : \ker B \to \ker(\sigma_1+1)
\end{align*}
and
$$S_{\sigma_1}(\lambda)=-
P_{\sigma_1}\circ C_1(\lambda)\circ I_{\sigma_1}: \ker
(\sigma_1+1) \to \ker (\sigma_1+1) \ \ .$$ For $\psi_1$ in
\eqref{express}, by Lemma \ref{l:inv}, there exists a unique
$\phi\in \ker(\sigma_1+1)$ such that $\psi_1=\phi-iG\phi$. As in
the derivation of \eqref{cc1}, we compare $\Psi_R$ with
$E(\phi,\la)$ and using Proposition \ref{p:nts1} we obtain
\begin{equation}\label{e:sm1}
\|C_1(\lambda) \psi_1 -\psi_2\| \le e^{-c R} \ \, .
\end{equation}
By the boundary condition in \eqref{e:boundary}, we have
\begin{equation*}
e^{-2i\lambda R}P_{\sigma_1}(\psi_1)=-P_{\sigma_1}(\psi_2) \ \, .
\end{equation*}
Combining this equation and (\ref{e:sm1}), we derive
$$
\|e^{2i\lambda R} S_{\sigma_1}(\lambda)\phi  -\phi\| \le e^{-cR}
$$
for $\phi \in \ker(\sigma_1 +1)$. We also define
\begin{align*}
I_{\sigma_2} &= I_+|_{\ker(\sigma_2+1)} : \ker(\sigma_2+1) \to
(\ker B)_- \ \ ,\\
P_{\sigma_2}&=\frac 12(1- \sigma_2) : \ker B \to \ker(\sigma_2+1)
\end{align*}
and
$$S_{\sigma_2}(\lambda):=-P_{\sigma_2} \circ C_2(\lambda)\circ I_{\sigma_2} :
\ker (\sigma_2+1) \to \ker (\sigma_2+1) \ \ ,$$ where $C_2(\la)$
is the scattering matrix defined from the generalized eigensection
over $M_{2,\infty}$. By the same way as above we can derive
$$
\|e^{2i\lambda R} S_{\sigma_2}(\lambda)\phi  -\phi\| \le e^{-cR}
$$
for $\phi \in \ker(\sigma_2 +1)$ . Now we introduce
\[
\Omega_i(R):=\{\, \rho\in\mathbb{R}-\{0\} \ | \ \det(e^{2i\rho
R}S_{\sigma_i}(\rho)-\Id)=0 \ ,\ |\rho| \le R^{-\kappa}\, \} \ \,
\]
for $i=1,2$. We repeat the argument used in \cite{Mu94} to prove
the corresponding result for \emph{s-values} of
$(\Dd_{i,R})_{P_i}$ noting all the argument for the involution
$C_i(0)$ in \cite{Mu94} holds for the involution $\sigma_i$, and
we obtain

\begin{thm}\label{t:small eigen2}For $i=1,2$,
let $\la_1(R)\le \la_2(R)\le \ldots \le \la_{p(R)}(R)$ be the
nonzero eigenvalues, counted to multiplicity, of
$(\Dd_{i,R})_{P_i}$ satisfying $|\la_k(R)|\le R^{-\kappa}$, and
let $\rho_1(R)\le \rho_2(R)\le \ldots \le \rho_{m(R)}(R)$ be the
nonzero element , counted to multiplicity, of $\Omega_i(R)$. Then
there exist $R_0$ and $c>0$, independent of $R$,  such that for
$R\ge R_0$, $p(R)=m(R)$ and
\[
|\la_k(R) -\rho_k(R)| \le e^{-cR} \quad \text{for}\quad
k=1,\ldots, p(R)\ \, .
\]
\end{thm}

 We now have the following
proposition

\begin{prop}\label{p:si_i}
There is a natural isomorphism
$$
\ker(\Dd_{i,R})_{P_i} \cong \ker_{L^2} \Dd_{i,\infty}\oplus
\ker(\sigma_i - 1)\cap L_i
$$
for $i=1,2$.
\end{prop}

\begin{proof} Let $\Psi\in \ker(\Dd_{1,R})_{P_1}$. Then the
section $\Psi$ satisfies $G(\partial_v+B)\Psi=0$ on the cylinder
$[0,R]_v\times Y$, and it has the following representation when
restricted to this cylinder
$$
\Psi= \phi_0 + \sum_{\mu_j >0} c_je^{-\mu_j v}\phi_j$$ where
$(\sigma_i-1)(\phi_0)=0$ . We use this expansion to extend $\Psi$
to a smooth section $\tilde{\Psi}$ on $M_{1,\infty}$ satisfying
$\Dd_{1,\infty}\tilde{\Psi}=0$. This means that $\tilde{\Psi}$
belongs to the space of the extended $L^2$-solutions of
$\Dd_{1,\infty}$. Hence $\phi_0$ is an element of $L_1$ . Let
$E(\phi_0,\lambda)$ be the generalized eigensection attached to
$\phi_0$. Then $\tilde{\Psi}-\frac 12 E(\phi_0,0)$ is square
integrable and $\Dd_{1,\infty}(\tilde{\Psi}-\frac12E(\phi_0,0))=0$
, and the map
$$
\Psi\longrightarrow (\tilde{\Psi}-\frac12E(\phi_0,0),\phi_0)$$
gives the expected isomorphism.
\end{proof}

The restriction of $S_{\sigma_i}:=S_{\sigma_i}(0)$ to
${\ker(\sigma_i+1) \cap \ker(C_i(0)+1)}$ is equal to the identity
map and
$$
\dim (\ker (\sigma_i+1) \cap \ker(C_i(0)+1))= \dim(\ker
(\sigma_i-1) \cap \ker(C_i(0)-1)) \ \,  .$$ It follows from
Proposition \ref{p:si_i} that the number of $(+1)$--eigenspace of
$S_{\sigma_i}:=S_{\sigma_i}(0)$ is equal to the dimension of the
subspace of $\ker(\Dd_{i,R})_{P_i}$ complementary to the subspace
$\ker_{L^2}\Dd_{i,\infty}$ for $i=1,2$ .

Now we define our model operator. Let $C : W \to W$ denote a
unitary operator acting on a $d$-dimensional vector space $W$ with
eigenvalues $e^{i\alpha_j}$ for $j=1,\ldots, d$. We introduce the
operator ${D}(C)$,
\begin{equation}\label{modeloperator}
D(C):=-i\frac12\frac{d}{du}:C^{\infty}(\mathbb{S}^1, E_C) \to
C^{\infty}(\mathbb{S}^1, E_C)
\end{equation}
where $E_C$ is the flat vector bundle over
$\mathbb{S}^1=\mathbb{R}/\mathbb{Z}$ defined by the holonomy
$\overline{C}$ , the complex conjugate of $C$ . The spectrum of
$D(C)$ is equal to\begin{equation}\label{e:spec} \{ \ \pi
k-\frac{1}{2}\alpha_j \ | \ k\in\mathbb{Z}, j=1,\ldots, d \ \} \
\, .
\end{equation}
The operators $C_{12}$ , $S_{\sigma_1}$ and $S_{\sigma_2}$ are the
unitary operators acting on finite dimensional vector spaces.
Hence we can define self-adjoint, elliptic operators
$D(C_{12}),D(S_{\sigma_1}), D(S_{\sigma_2})$ on $\mathbb{S}^1$.

\begin{thm}\label{t:model}
Assume that all the \emph{e-values} of $\Dd_R$ are zero
eigenvalues. Let $\la_1(R)\le \la_2(R)\le \ldots \le
\la_{p(R)}(R)$ be the nonzero eigenvalues, counted to
multiplicity, of $\Dd_R$ satisfying $|\la_k(R)|\le R^{-\kappa}$,
and let $\la_1\le \la_2\le \ldots \le \la_{n(R)}$ be the nonzero
eigenvalues, counted to multiplicity, of $D(C_{12})$ satisfying
$|\la_k|\le 2R^{1-\kappa}$. Then there exist $R_0$ and $C>0$,
independent of $R$,  such that for $R\ge R_0$, $p(R)=n(R)$ and
\[
|2R\la_k(R) -\la_k| \le C\, R^{-\kappa} \quad \text{for}\quad
k=1,\ldots, p(R)\ \, .
\]
The similar statement holds for $(\Dd_{i,R})_{P_i}$ and
$D(S_{\sigma_i})$ with the relation
\[
|R\la_k(R) -\la_k| \le C\, R^{-\kappa} \quad \text{for}\quad
k=1,\ldots, p_i(R)\ \,
\]
where $p_i(R)$ is the number of \emph{s-values} of
$(\Dd_{i,R})_{P_i}$ with $|\la_k(R)| \le R^{-\kappa}$.
\end{thm}

\begin{proof}
First we introduce
\[
\Omega^*(R):=\{\, \rho\in\mathbb{R}-\{0\} \ | \ \det(e^{4i\rho
R}C_{12}-\Id)=0 \ ,\ |\rho| \le R^{-\kappa}\, \} \ \, .
\]
By definition, this set consists of the nonzero solution
$\rho^*_{j,k}$ of
\begin{equation}\label{sol*}
4\la R +\alpha_j(0)=2\pi k \quad \text{with} \quad |\la| \le
R^{-\kappa} \ \,
\end{equation}
where $e^{i\alpha_j(0)}$ for $j=1,\ldots, \frac{h_Y}2$ are the
eigenvalues of $C_{12}=C_{12}(0)$. Now, for an element
$\rho_{j,k}$ in $\Omega(R)$ defined in \eqref{OmegaR},  one can
show (as near \eqref{e:sol}) that if $R\gg0$ there is the
corresponding solution $\rho_{j,k}^*$ of
\[
4\la R +\alpha_j(0)=2\pi k \quad \text{with} \quad |\la| \le
R^{-\kappa}+R^{-1-\kappa} \ \, ,
\]
noting $|\alpha_j(\la)-\alpha_{j}(0)| \le cR^{-\kappa}$ for a
positive constant $c$. Since $|\rho_{j,k}-\rho^*_{j,k}| \le
cR^{-1-\kappa}$, this gives a one to one correspondence from
$\Omega(R)$ to $\Omega^*(R_0)$ with
$R_0^{-\kappa}=R^{-\kappa}+R^{-1-\kappa}$ for $R\gg 0$. Now, let
us observe that for any pair of $\rho^*_{j,k}\neq \rho^*_{j',k'}$
in $\Omega^*(R)$, $|\rho^*_{j,k}-\rho^*_{j',k'}| \ge {a_0}R^{-1} $
for a positive constant $a_0$. Hence, for $R\gg0$, this implies
that $\Omega^*(R)=\Omega^*(R_0)$ with
$R_0^{-\kappa}=R^{-\kappa}+R^{-1-\kappa}$. In conclusion, there is
a one to one correspondence between $\Omega(R)$ and $\Omega^*(R)$
for $R\gg 0$ with the relation
\[
|\rho_k -\rho_k^* | \le c R^{-1-\kappa}
\]
where $\rho_1\le \ldots \le \rho_{m(R)}$ ($\rho^*_1\le \ldots \le
\rho^*_{n(R)}$) denotes the elements, counted to multiplicity, of
$\Omega(R)$ ($\Omega^*(R)$). For $\rho^*\in\Omega^*(R)$, the map
$\rho^*\to 2R\rho^*$ gives a one to one correspondence from
$\Omega^*(R)$ to the subset of the eigenvalues $\la_k$ of
$D(C_{12})$ with $|\la_k|\le 2R^{1-\kappa}$. Now, applying Theorem
\ref{t:small eigen1} completes the proof for \emph{s-values} for
$\Dd_R$. The case of $(\Dd_{i,R})_{P_i}$ can be proved in the same
way.

\end{proof}

\section{Large time contribution}\label{four}

In this section, we prove the following proposition

\begin{prop}\label{p:lc}
\begin{multline*}\label{e:lc}
\lim_{R\to\infty} \int^{\infty}_{R^{2-\e}}
t^{-1}[\Tr(e^{-t\Dd_R^2} - e^{-t(\Dd_{1,R})_{P_1}^2} -
e^{-t(\Dd_{2,R})_{P_2}^2}) -h ]\ dt - h(\gamma - \e{\cdot}\log
R)\\
= \frac{d}{ds}\Big|_{s=0} \frac{1}{\Gamma(s)}\int^{\infty}_0
t^{s-1} [\Tr(e^{-t\frac{1}{4}D(C_{12})^2}-e^{-tD(S_{\sigma_1})^2}-
e^{-tD(S_{\sigma_2})^2})-h] \ dt
\end{multline*}
where $h=\dim(L_1\cap L_2) - \dim(L_1\cap \ker(\sigma_1-1)) -
\dim(L_2\cap \ker(\sigma_2-1))$ .
\end{prop}

Recalling
\[
\frac{d}{ds}\Big|_{s=0} \z^R_l(s)= \int^{\infty}_{R^{2-\e}}t^{-1}
[\, \Tr(e^{-t\Dd_R^2} - e^{-t(\Dd_{1,R})_{P_1}^2} -
e^{-t(\Dd_{2,R})_{P_2}^2})- h\, ] \ dt \ \, ,
\]
Proposition \ref{p:lc} immediately implies that the large time
contribution to the adiabatic decomposition formula for the
$\z$-determinant is equal to
$$\frac{\zd\frac 14 D(C_{12})^2}{\zd D(S_{\sigma_1})^2
{\cdot}\zd D(S_{\sigma_2})^2} \ \ .$$

We start with the following result,

\begin{prop}\label{p:s1}
The following equality holds,
\begin{multline*}
\lim_{R \to \infty}\Bigr(\frac{d}{ds}\Big|_{s=0}
\frac{1}{\Gamma(s)}
\int^{R^{-\e}}_0t^{s-1}[ \Tr(e^{-t\frac{1}{4}D(C_{12})^2}\\
\quad -e^{-tD(S_{\sigma_1})^2} -e^{-tD(S_{\sigma_2})^2})-h]\, dt
+h(\gamma - \e{\cdot}\log R)\Bigr) = 0 \ \ .
\end{multline*}

\end{prop}

\begin{proof}
Recalling the definition of $D(C)$ in \eqref{modeloperator}, we
can see that if $\Ll$ is one of $D(C_{12})^2$ ,
$D(S_{\sigma_1})^2$, $D(S_{\sigma_1})^2$, then
$$
\Tr(e^{-t\Ll}) \sim \sqrt{\frac {\pi}t}\, \frac{h_Y}{2} +
O(\sqrt{t}) \quad \text{near}\quad t=0 \ \
$$
since $h_Y=2\dim \Gg_-=2\dim \ker(\sigma_i+1)$. Hence there exists
a constant $c_1$ such that
\begin{equation}\label{e:sd1}
|\Tr(e^{-t\frac{1}{4}D(C_{12})^2}
-e^{-tD(S_{\sigma_1})^2}-e^{-tD(S_{\sigma_2})^2})| <
{c_1}{\sqrt{t}} \quad \text{near}\quad t=0 \ \, .
\end{equation}
This allows us to estimate
\begin{multline*}
\Big|\, \frac{d}{ds}\Big|_{s=0}\frac{1}{\Gamma(s)}\int^{R^{-\e}}_0
t^{s-1} \Tr(e^{-t\frac{1}{4}D(C_{12})^2} -
e^{-tD(S_{\sigma_1})^2}-
e^{-tD(S_{\sigma_2})^2})\ dt\, \Big| \\
\le c_2{\cdot}\int_0^{R^{-\e}}\frac {dt}{\sqrt{t}}
=2c_2{\cdot}R^{-\frac {\e}2} \ \ .
\end{multline*}
\bigskip
Combining this with
$$
\frac{d}{ds}\Big|_{s=0} \frac{h}{\Gamma(s)} \int^{R^{-\e}}_{0}
t^{s-1}dt  = h(\gamma - \e \cdot \log R) \ \ $$ completes the
proof.
\end{proof}

 It follows from Proposition \ref{p:s1} that Proposition
\ref{p:lc} is equivalent to the following equation
\begin{multline}\label{e:claim1}
\lim_{R\to\infty} \int^{\infty}_{R^{2-\e}} t^{-1}[\,
\Tr(e^{-t\Dd_R^2} - e^{-t(\Dd_{1,R})_{P_1}^2} -
e^{-t(\Dd_{2,R})_{P_2}^2}) - h\, ] \ dt \\
=\lim_{R\to\infty}\frac{d}{ds}\Big|_{s=0} \frac{1}{\Gamma(s)}
\int^{\infty}_{R^{-\e}}t^{s-1} [\,
\Tr(e^{-t\frac{1}{4}D(C_{12})^2}
-e^{-tD(S_{\sigma_1})^2}-e^{-tD(S_{\sigma_2})^2})-h\, ]\ dt.
\end{multline}
Now using change variables we obtain
\begin{align*}
&\int^{\infty}_{R^{2-\e}} t^{-1}[\, \Tr(e^{-t\Dd_R^2} -
e^{-t(\Dd_{1,R})_{P_1}^2} - e^{-t(\Dd_{2,R})_{P_2}^2})- h\, ]\ dt
\\
=&\int^{\infty}_{R^{-\e}}t^{-1} [\, \Tr(e^{-tR^2\Dd_R^2} -
e^{-tR^2(\Dd_{1,R})_{P_1}^2} - e^{-tR^2(\Dd_{2,R})_{P_2}^2})- h \,
]\ dt \ \ .
\end{align*}
Then the equality (\ref{e:claim1}) is equivalent to
\begin{multline}\label{e:lc1}
\lim_{R\to\infty}\int^{\infty}_{R^{-\e}} t^{-1}
[\,\Tr(e^{-tR^2\Dd_R^2} - e^{-tR^2(\Dd_{1,R})_{P_1}^2} -
e^{-tR^2(\Dd_{2,R})_{P_2}^2})- h\,] \ dt\\
=\lim_{R\to\infty}\frac{d}{ds}\Big|_{s=0} \frac{1}{\Gamma(s)}
\int^{\infty}_{R^{-\e}}t^{s-1} [\,
\Tr(e^{-t\frac{1}{4}D(C_{12})^2}
-e^{-tD(S_{\sigma_1})^2}-e^{-tD(S_{\sigma_2})^2})-h\, ]\ dt.
\end{multline}

Now we split
$$\Tr(e^{-tR^2\Dd_R^2} - e^{-tR^2(\Dd_{1,R})_{P_1}^2} -
e^{-tR^2(\Dd_{2,R})_{P_2}^2}) - h $$ into two parts
\begin{align*}
&\Tr_{I,R}(e^{-tR^2\Dd_R^2} - e^{-tR^2(\Dd_{1,R})_{P_1}^2} -
e^{-tR^2(\Dd_{2,R})_{P_2}^2})\ \ ,\\
&\Tr_{II,R}(e^{-tR^2\Dd_R^2} - e^{-tR^2(\Dd_{1,R})_{P_1}^2} -
e^{-tR^2(\Dd_{2,R})_{P_2}^2})
\end{align*}
where $\Tr_{I,R}(\cdot)$ ($\Tr_{II,R}(\cdot)$) is  the part of the
trace restricted to the eigenvalues of $R^2\Dd_{R}^2$ ,
$R^2(\Dd_{1,R})_{P_1}^2$ , $R^2(\Dd_{2,R})_{P_2}^2$ which are
larger (smaller or equal to) than $R^{\frac 12}$.  The next
proposition shows that $\Tr_{I,R}(\cdot)$ can be neglected as
$R\to\infty$,

\bigskip

\begin{prop}\label{p:int}
We have the following estimate
$$\int^{\infty}_{R^{-\e}} t^{-1}  \Tr_{I,R}(e^{-tR^2\Dd_R^2} -
e^{-tR^2(\Dd_{1,R})_{P_1}^2} - e^{-tR^2(\Dd_{2,R})_{P_2}^2})\ dt
\le {c_1}e^{-c_2 R^{\frac 12-\e}}$$ for some positive constants
$c_1,c_2$.
\end{prop}

\begin{proof}
Let $\la_{k_0}(R)^2$ denote the smallest \emph{large} eigenvalue
of $\Dd_R^2$, that is, the smallest one with $\la_{k_0}(R)^2 >
R^{-\frac32}$. We now estimate $\Tr_{I,R}(e^{-tR^2\Dd_{R}^2})$ as
follows;
\begin{multline*}
\Tr_{I,R} (e^{-tR^2\Dd_{R}^2}) = \sum_{\la_k^2 >
R^{-\frac32}}e^{-tR^2\la_k^2} = \sum_{\la_k^2 > R^{-\frac32}} e^{-(tR^2-1)\la_k^2}e^{-\la_k^2} \\
\le e^{-(tR^2-1)\la_{k_0}^2}\sum_{\la_k^2 >
R^{-\frac32}}e^{-\la_k^2} \le e^{-(tR^2-1)\la_{k_0}^2}\Tr (
e^{-\Dd_{R}^2})\\
\le  b_1Re^{-(tR^2-1)R^{-\frac 32}} \le b_1Re^{-b_2t\sqrt{R}} \ \
,
\end{multline*}
for some positive constants $b_1$, $b_2$ . Hence we have
\begin{multline*}
\int^{\infty}_{R^{-\e}}{t^{-1}}\Tr_{I,R}(e^{-tR^2\Dd_{R}^2}) \ dt
\le \int^{\infty}_{R^{-\e}}{{t^{-1}}}b_1Re^{-b_2t\sqrt{R}}\
dt\\
 \le {\frac {b_1}{b_2}} R^{{\frac 12} +\e} \int_{b_2R^{\frac 12 -
\e}}^{\infty}e^{-v}dv \le b_3e^{-b_4R^{\frac 12 - \e}}  \ \ .
\end{multline*}
The trace $\Tr_{I,R}(e^{-tR^2(\Dd_{i,R})_{P_i}^2})$ for $i=1,2$
can be estimated in the same way. This completes the proof.

\end{proof}

\bigskip

\noi We also split
$\Tr(e^{-t\frac{1}{4}D(C_{12})^2}-e^{-tD(S_{\sigma_1})^2}-
e^{-tD(S_{\sigma_2})^2})-h$ into two parts
\begin{align*}
&\Tr_{I,R}(e^{-t\frac{1}{4}D(C_{12})^2}
-e^{-tD(S_{\sigma_1})^2}-e^{-tD(S_{\sigma_2})^2}) \ \ ,\\
&\Tr_{II, R}(e^{-t\frac{1}{4}D(C_{12})^2} -
e^{-tD(S_{\sigma_1})^2}-e^{-tD(S_{\sigma_2})^2})
\end{align*}
where $\Tr_{I,R}(\cdot)$ ($\Tr_{II,R}(\cdot)$) is taken over the
nonzero eigenvalues of $\frac 14 D(C_{12})^2$,
$D(S_{\sigma_1})^2$, $D(S_{\sigma_2})^2$ which are larger (smaller
or equal to) than $R^{\frac 12}$. The following proposition
corresponds to Proposition \ref{p:int} and its proof is
essentially the same as the proof of Proposition \ref{p:int}.

\begin{prop}\label{p:int2}
There exist positive constants $c_{1},c_{2}$ such that
$$
\int^{\infty}_{R^{-\e}}
t^{-1}\Tr_{I,R}(e^{-t\frac{1}{4}D(C_{12})^2}
-e^{-tD(S_{\sigma_1})^2}- e^{-tD(S_{\sigma_2})^2}) \le
c_{1}e^{-c_{2}R^{\frac 12-\e}} \ \ .$$
\end{prop}

By Propositions \ref{p:int} and \ref{p:int2}, we can see that the
equality (\ref{e:lc1}) is equivalent to
\begin{multline}\label{e:lc2}
 \lim_{R\to\infty} \biggl(\int^{\infty}_{R^{-\e}}t^{-1}
  \Tr_{II,R}(e^{-tR^2\Dd_R^2} - e^{-tR^2(\Dd_{1,R})_{P_1}^2}
- e^{-tR^2(\Dd_{2,R})_{P_2}^2}) \ dt \\
 \qquad\qquad\qquad -\int^{\infty}_{R^{-\e}}t^{-1} \Tr_{II,R}
(e^{-t\frac{1}{4}D(C_{12})^2} -e^{-tD(S_{\sigma_1})^2}-
e^{-tD(S_{\sigma_2})^2}) \ dt \biggr) \ =\  0 \ \ .
\end{multline}
The equation (\ref{e:lc2}) is a consequence of the next result

\begin{prop}\label{p:diff}
For sufficiently large $R$, there exist positive constants
$c_{1},c_{2}$ independent of $R$ and $t$, such that
\begin{align*}
&|\ \Tr_{II,R}(e^{-tR^2\Dd_R^2})
-\Tr_{II,R}(e^{-t\frac{1}{4}D(C_{12})^2})\ | \ \le \
c_{1}{t}{R^{-\frac 14}} e^{-c_{2} t} \ \ ,\\
&|\ \Tr_{II,R}(e^{-tR^2(\Dd_{i,R})_{P_i}^2})-\Tr_{II,R}
(e^{-tD(S_{\sigma_i})^2})\ |\ \le \  c_{1}{t}{R^{-\frac 14}}
e^{-c_{2}t}
\end{align*}
for any $t>0$.
\end{prop}

\begin{proof}
We use the analysis of \emph{s-values} developed in Section
\ref{three}. We fix $\kappa = \frac 34$ . It follows from Theorem
\ref{t:model} that, for any eigenvalue $\lambda(R)$ of $\Dd_R$
with $|\lambda(R)|\le R^{-\frac 34}$ there exists an analytic
function $\alpha(\lambda)$ such that
$$
R\lambda(R)=\lambda_j +\frac14\lambda(R)\alpha(\lambda(R)) +
O(e^{-cR})
$$
where $\lambda_j$ is an eigenvalue of $\frac 12 D(C_{12})$ with
$|\lambda_j| \le R^{\frac 14}$. Therefore, there exist functions
$c(R), d(R)$ and a constant $C>0$ such that
$$R^2\lambda(R)^2= \lambda_j^2 + \la_j\frac{c(R)}{R^{\frac 34}}
+ \frac{d(R)}{R^{\frac32}} \quad\text{with}\quad |c(R)| \le C\ ,
\quad |d(R)|\le C \ \ ,$$ for any sufficiently large $R$ . We use
the elementary inequality $|e^{-\lambda}-1|\le
|\lambda|e^{|\lambda|}$ to get
\begin{multline*}
|e^{-tR^2\lambda(R)^2}-e^{-t\lambda_j^2}|= |e^{-t\lambda_j^2}
(e^{-t[R^2\lambda(R)^2-\lambda_j^2]}-1)|\\
\le \Big(\frac{|\lambda_j c(R)|}{R^{\frac 34}}+
\frac{|d(R)|}{R^{\frac32}}\Big)\,t\, e^{-(\lambda_j^2-\frac{|
\lambda_j c(R)|}{R^{\frac34}}-\frac{|d(R)|}{R^{\frac32}})t} \le
\frac{ C}{R^{\frac 12}}\, t\, e^{-\frac 12 \lambda_j^2t} \ \
\end{multline*}
for $R\gg 0$. In the last inequality we used the fact that
$|\lambda_j| \le R^{\frac 14}$ . Let us fix a sufficiently large
$R$. We take the sum over finitely many eigenvalues $\lambda_j^2$
of $\frac 14 D(C_{12})^2$ with $\lambda_j^2 \le R^{\frac 12}$, and
obtain
$$|\Tr_{II,R}(e^{-tR^2\Dd_R^2}) -
\Tr_{II,R}(e^{-t\frac{1}{4}D(C_{12})^2})|\le C \frac{t}{R^{\frac
12}} \sum_{\lambda_j^2\le R^{\frac 12}} e^{-\frac 12\lambda_j^2t}
\ \ .$$  The operator $\frac14 D(C_{12})^2$ is a Laplace type
operator over $S^1$ , hence  the number of eigenvalues
$\lambda_j^2$ of $\frac14 D(C_{12})^2$ with $\lambda_j^2\le
R^{\frac 12}$ can be estimated by $R^{\frac 14}$. Therefore, we
have
$$C \frac{t}{R^{\frac 12}} \sum_{\lambda_j^2 \le R^{\frac 12}}
e^{-\frac 12\lambda_j^2t} \le c_1\frac{t}{R^{\frac 12}}R^{\frac
14} e^{-\frac 12 \lambda^2_1 t} \ \,$$ where $\la_1^2$ denotes the
first nonzero eigenvalue of $\frac 14 D(C_{12})^2$ . Note that
$c_1$ and $\lambda_1^2$ are independent of $R$.  This proves the
first claim putting $c_2=\frac12\la_1^2$. The proof of the second
claim goes in the same way.
\end{proof}

The proof of Proposition \ref{p:lc} is now complete.

{\flushleft \textbf{Proof of Theorem \ref{t:apssplit}:}} Now
Proposition {\ref{p:s}} and Proposition \ref{p:lc} give us the
following equality,
\begin{align*}
 &\lim_{R \to \infty}\Bigr({\z_s^R}'(0) + h(\g + (2-
\e){\cdot}\log R)
+ {\z_l^R}'(0) - h(\g - \e{\cdot}\log R)\Bigr)\\
= & \frac{d}{ds}\Big|_{s=0} \frac{1}{\Gamma(s)}\int^{\infty}_0
t^{s-1} [\,
\Tr(e^{-t\frac{1}{4}D(C_{12})^2}-e^{-tD(S_{\sigma_1})^2}-
e^{-tD(S_{\sigma_2})^2})-h \, ]\ dt \notag  \\
 & + \z_{B^2}(0)\cdot  \log 2 \ \ . \notag
\end{align*}
By an elementary computation (for instance, see Proposition 2.2 in
\cite{LP3}), we can derive
\begin{align}\label{e:sindet}
\zd \frac14D(C_{12})^2 \ =& \ 2^{h_Y+2h_M}\, {\det}^* \Big( \frac
{2\,\Id - C_{12} - C_{12}^{-1}}4 \Big)\ \, , \notag \\ \zd
D(S_{\sigma_i})^2 \ =& \ 2^{h_Y}\, {\det}^* \Big( \frac {2\, \Id -
S_{\sigma_i} - S_{\sigma_i}^{-1}}4 \Big)\ \, .
\end{align}
Combining all these equalities provides us with the final formula
(\ref{e:main}) in Theorem \ref{t:apssplit}.

\section{A proof of the decomposition formula of the $\eta$-invariant}
\label{five}

In this section we offer a new proof of the decomposition formula
for the $\eta$-invariant. This formula  has been proved by several
authors (see \cite{Bu95}, \cite{Dfr94}, \cite{HMaM95},
\cite{Mu96}, \cite{KPW95}, \cite{BL98}, \cite{KL01}, \cite{LP1})
and the proof we discuss in this section is not the simplest one.
Still we believe that it is worthy to present the \emph{scattering
approach} to the decomposition of the $\eta$-invariant. The key in
our proof is to show that the {\it scattering data} provides us
with the contribution given by the boundary conditions in the
decomposition formula for the $\eta$-invariant.

Let us remind the reader that the $\eta$-function of a Dirac
operator $\Dd$ on a closed manifold $M$, introduced in
\cite{AtPaSi75}, is defined as
$$\eta_\Dd(s) = \sum_{\la_k \ne 0}\sign(\la_k)|\la_k|^{-s} \ \ ,$$
where the sum is taken over all nonzero eigenvalues of $\Dd$ . The
$\eta$-function is well-defined for $\Re(s)$ large and it has a
meromorphic extension to the whole complex plane and $s=0$ is a
regular point, hence $\eta_{\Dd}(0)$ is well-defined. Following
\cite{AtPaSi75} we introduce the $\eta$-invariant of $\Dd$ as
\begin{equation}\label{e:eta1}
\eta(\Dd) = \frac 12{\cdot}(\eta_{\Dd}(0) + \dim \ker  \Dd) \ \, .
\end{equation}

Now, let us assume that we have a decomposition of a closed
odd-dimensional manifold $M$ to $M_1\cup M_2$ in the way described
in the introduction. For $\Dd_i:=\Dd|_{M_i}$, we impose the
boundary conditions given by the generalized APS spectral
projections $P_i$ defined in \eqref{e:id1}. Then the
$\eta$-function of $(\Dd_{i})_{P_i}$ is also well-defined and it
has the same properties as the $\eta$-function of the Dirac
operator on a closed manifold, in particular, the $\eta$-function
of $(\Dd_i)_{P_i}$ is regular at $s=0$. Hence, we can define the
$\eta$-invariant of $(\Dd_i)_{P_i}$ as in \eqref{e:eta1}. The
following result was proved by several authors as we remarked
above,

\begin{thm}\label{t:eta-1}
The following formula holds,
\begin{equation}\label{e:eta2}
\eta(\Dd) = \eta((\Dd_1)_{P_1}) + \eta((\Dd_2)_{P_2}) +
\eta(\Dd;\s_1,\s_2) \ \mod  \mathbb{Z} \ \, ,
\end{equation}
where $\eta(\Dd;\s_1,\s_2)$ denotes the $\eta$-invariant of the
operator $\Dd=G(\partial_u + B)$ over $N\cong [-1,1]\times Y$,
subject to the boundary condition $P_2$ at $u=-1$ and $P_1$ at
$u=1$.
\end{thm}

For the involution $\sigma_i$ which defines $P_i$ in
\eqref{e:id1}, let us observe that
$$U = \s_1\s_2 : \ker B \to \ker B$$
is the unitary operator, such that $UG = GU$, $\det U = 1$ and
$U^*=\sigma_1U\sigma_1$. It follows that the spectrum of $U$ is
invariant under complex conjugation. Moreover, the maps $U_{\pm} =
U|_{\Gg_{\pm}} : \Gg_{\pm} \to \Gg_{\pm}$ are well-defined. The
following result proved in Section 2 of \cite{LeWo96} was the key
ingredient in the proof of Theorem \ref{t:eta-1}.

\vskip 5mm

\begin{prop}\label{t:eta1}
We have the following formulas
\begin{equation}\label{e:eta4}
\eta(\Dd;\s_1,\s_2) = -{\frac {1}{2\pi i}\log \det(-U_+) } \ \mod
 \mathbb{Z} \ \,  .
\end{equation}

\end{prop}

 One way to prove the decomposition formula (\ref{e:eta2}) is to
use the \emph{adiabatic analysis} we developed in the proof of
Theorem \ref{t:apssplit}. This analysis easily gives us the
following theorem,

\begin{thm}\label{t:eta2}
The following formula for the $\eta$-invariant holds,
\begin{multline}\label{e:eta5}
 \eta(\Dd) - \eta((\Dd_{1})_{P_1}) - \eta((\Dd_{2})_{P_2})\\
 =  \eta(D(C_{12}))-\eta(D(S_{\s_1}))-\eta(D(S_{\s_2}))
\ \mod   \mathbb{Z}\ \, . \end{multline}
\end{thm}

\begin{proof}
We repeat the corresponding argument to derive Theorem
\ref{t:apssplit} for the $\eta$-invariant to obtain the expected
formula
\begin{multline*}
\lim_{R \to \infty}\{\eta(\Dd_R) - \eta((\Dd_{1,R})_{P_1}) -
\eta((\Dd_{2,R})_{P_2})\}\\ =
\eta(D(C_{12}))-\eta(D(S_{\s_1}))-\eta(D(S_{\s_2})) \ \mod
\mathbb{Z}.
\end{multline*}
Now, we use the fact that $\eta(\Dd_R), \eta((\Dd_{i,R})_{P_i})$
are independent of $R$ modulo integer (see Proposition 2.16 of
\cite{Mu94}) to complete the proof.
\end{proof}

Now we need to show
$$\eta(\Dd;\s_1,\s_2)=\eta(D(C_{12}))-\eta(D(S_{\s_1}))-\eta(D(S_{\s_2}))
 \ \mod
\mathbb{Z}\ \, . $$ For this, we observe the followings: The
scattering matrix $C_i = C_i(0)$ can be represented in the
following way,
$$C_i =
\left(\begin{matrix} 0 & C(i)_- \\
C(i)_+ & 0 \end{matrix} \right ) \ \  \text{where} \ \
C(i)_{\pm}C(i)_{\mp} = \Id \ \ ,
$$
with respect to the decomposition  $\ker B = \Gg_+ \oplus \Gg_-$.
We see that
\begin{equation}\label{c1}
C_{12} = C(1)_+C(2)_- : \Gg_- \to \Gg_- \ \ . \end{equation}
Similar formulas hold for the involutions $\s_i$ and we have
\begin{align*}
S_{\s_1}=&-P_{\sigma_1}\circ C_1\circ I_{\sigma_1} = -{\frac 12}
\left(\begin{matrix} \Id & -\s(1)_- \\
-\s(1)_+ & \Id \end{matrix}\right)
\left(\begin{matrix} 0 & C(1)_- \\
C(1)_+ & 0 \end{matrix}\right)
\left(\begin{matrix} 2\mathrm{Id} & 0 \\
0 & 0 \end{matrix}\right)\\
 =&\left(\begin{matrix} \s(1)_-C(1)_+ & 0 \\
-C(1)_+ & 0 \end{matrix}\right) \ \ .
\end{align*}
We can also see that  each element of $\ker(\s_1 + 1)$ is
represented in the form $\left(\begin{matrix} f\\
-\s(1)_+f\end{matrix}\right)$ for some $f \in \Gg_+$ . This allows
us to represent the map $S_{\s_1}$ over $\ker(\s_1 + 1)$ as
\begin{align*}
S_{\s_1}\left(\begin{matrix} f\\
                      -\s(1)_+f\end{matrix}\right) =
\left(\begin{matrix} \s(1)_-C(1)_+ & 0 \\
-C(1)_+ & 0 \end{matrix}\right)\left(\begin{matrix} f\\
                      -\s(1)_+f\end{matrix}\right)
                      =
\left(\begin{matrix} \s(1)_-C(1)_+f\\
                     -\s(1)_+\s(1)_-C(1)_+f\end{matrix}\right)
\  \, .
\end{align*}
Therefore, from the spectral point of view, the operator
$S_{\s_1}$ is equal to the operator
$$\s(1)_-C(1)_+ : \Gg_+ \to \Gg_+ \ \ ,$$
or equivalently to the operator
$$C(1)_+\s(1)_- : \Gg_- \to \Gg_- \ \ .$$
The corresponding analysis for the operator $S_{\s_2}$ implies
that $S_{\s_2}$ is equivalent to
$$C(2)_-\s(2)_+ : \Gg_+ \to \Gg_+ \ \ \text{or} \ \
\s(2)_+C(2)_- : \Gg_- \to \Gg_- \ \ .$$ Combining \eqref{c1} and
these, we obtain
\begin{align}\label{det}
 \frac {\det(C_{12})}{\det(S_{\s_1})\det(S_{\s_2})} =
 \det  (\s(1)_+\s(2)_-) \ \ .
\end{align}
For the operator $D(C)$ on $S^1$ defined by a unitary map $C$ in
\eqref{modeloperator},
\begin{equation}\label{e:etas}
\eta(D(C)) = - \frac {1}{2\pi i}\log \det(-\overline{C}) \ \mod
\mathbb{Z} \ \, .
\end{equation}
(see Theorem 2.1 and Lemma 2.3 in \cite{LeWo96}). If we combine
\eqref{det} and \eqref{e:etas}, we have
\begin{align*}
\eta(D(C_{12}))-\eta(D(S_{\sigma_1}))-\eta(D(S_{\sigma_2})) =
-\frac{1}{2\pi i}{\log \det(-\overline{\s(1)_+\s(2)_-}})
  \ \mod  \mathbb{Z} \ \ .
\end{align*}
Noting $\det(-\overline{\s(1)_+\s(2)_-})=\det(-\s(1)_-\s(2)_+)$,
this and Proposition \ref{t:eta1} end the proof of the following
theorem,

\begin{thm}\label{t:eta3}
\begin{equation}\label{e:eta61}
\eta(\Dd;\sigma_1,\sigma_2)=\eta(D(C_{12}))-\eta(D(S_{\sigma_1}))-
\eta(D(S_{\sigma_2})) \ \mod  \mathbb{Z} \ \, .
\end{equation}
\end{thm}

\end{document}